\documentclass{birkjour}
\usepackage{amsmath, amssymb, amsthm}

\newcommand
{\Rd}{\mathrm{Dom}}

\newcommand
{\R}{{\mathbf R}}

\newcommand
{\T}{\mathrm{graph}}

\newcommand
{\I}{\mathrm{Im}}

\newcommand
{\E}{\mathrm{pr}}

\newcommand
{\C}{{\mathcal C}}

\newcommand
{\M}{{\mathcal M}}

\newcommand
{\D}{{\mathbf D}}

\newcommand
{\B}{{\mathcal H}}

\newcommand
{\G}{{\mathcal G}}

\newcommand
{\F}{{\mathcal F}}

\newcommand
{\K}{{\mathcal L}}

\newcommand
{\A}{{\mathcal H}}

\newcommand
{\J}{{\mathcal D}}

\newcommand
{\U}{{\mathcal U}}

\newcommand
{\X}{{\mathcal T}}

\newcommand
{\V}{{\mathcal V}}

\newcommand
{\Y}{{\mathcal S}}

\newcommand
{\W}{{\mathcal W}}

\newtheorem{theorem}{Theorem}[section]
\newtheorem{proposition}[theorem]{Proposition}
\newtheorem{lemma}[theorem]{Lemma}

\theoremstyle{definition}
\newtheorem{definition}[theorem]{Definition}
\theoremstyle{remark}
\newtheorem{remark}[theorem]{Remark}
\numberwithin{equation}{section}

\renewcommand{\theequation}{\thesection.\arabic{equation}}

\begin{document}
	
	\title[Extensions of Nonlinear Smooth Operators]{Locally Self-adjoint Extensions of\\ Nonlinear Smooth Operators and \\Abstract Boundary Conditions}
	
	\author[L. Zelenko]
	{ Leonid Zelenko} 
	
	%\date{}
	\address{%
		Department of Mathematics \\
		University of Haifa  \\
		Haifa 31905  \\
		Israel}
	\email{zelenko@math.haifa.ac.il}

\begin{abstract}
In a real Hilbert spaces $\B$ a smooth operator $F$ is studied,
whose derivative $F^\prime(x)$ at each point $x\in\Rd(F)$ is a
symmetric operator. In terms of abstract boundary conditions locally
self-adjoint extensions of this operator are described. We use some 
concepts and facts from symplectic differential geometry.
\end{abstract}



  \maketitle

\tableofcontents
 
\section{\bf Introduction}
\setcounter{equation}{0}

The role played by the theory of extensions of linear operators in 
Hilbert space in the analysis of linear boundary value problems 
is well known [D-Sh], [Mr]. Self-adjoint boundary value problems 
are very important from the physical point of view. They describe 
conservative systems, i.e. systems which are subordinated to some 
conservation law. For example, the boundary value problems for the 
Euler equations generated by a wide class of variation problems 
are self-adjoint in the above sense. The problem of description 
of linear self-adjoint boundary conditions 
for linear differential equations 
(ordinary or partial) is reduced to the problem of description 
of self-adjoint extensions of the corresponding linear symmetric 
operators in a Hilbert space. At first the theory of such 
extensions was constructed by John von Neumann \cite{VN} and afterwards 
it was developed in papers of K. Friedrichs, H. Freudenthal, 
J. Calkin, M. Krein and other authors (\cite{Frd}, \cite{Frdn}, \cite{Cl}, \cite{Kr}). 
Some generalizations of Friedrichs extension to the nonlinear case were constructed in \cite{Pet} and \cite{Zl-1}.

The self-adjoint extensions mentioned above can be described 
by means of so called "abstract boundary conditions". Conditions 
of such kind are given by the Calkin Theorem \cite{Cl}. Let us formulate 
it.

\textbf{Calkin Theorem}. Let $F_0$ be a closed linear symmetric operator 
in a Hilbert space $\A$ with a dense domain $\Rd(F_0)$ and with finite 
defects $\{m,m\}$. Denote by $\langle\cdot,\cdot\rangle$ the following bilinear form 
defined on $\Rd(F_0^\star)$:
\begin{equation*}
\langle v_1,v_2\rangle=(F_0^\star v_1,v_2)-(v_1,F_0^\star v_2).
\end{equation*}
A linear set $D\subseteq\A$ is the domain of a self-adjoint extension 
$F$ of the operator $F_0$ if and only if there exists such a system 
of elements $v_1,v_2,\dots,v_m\in\A$, that 

(a) $v_1,v_2,\dots,v_m$ are linearly independent modulo $\Rd(F_0)$;

(b) $\langle v_i,v_j\rangle=0\;\;\forall i,j\in\{1,2,\dots,m\}$;

(c) $D=\{f\in\Rd(F_0^\star):\;\langle f,v_i\rangle=0         
\;\;\forall i\in\{1,2,\dots,m\}\}$

The complete description of self-adjoint extensions of 
linear symmetric ordinary 
differential operators (regular and singular) on the base of the 
Calkin Theorem was carried out in the book of M. A. Naimark [Nm]
(see also \cite{D-Sh}).
 
 Consider the 
 orthogonal sum:
 \begin{equation}\label{b7}
 \B^2 = \B\oplus \B.
 \end{equation}
 We denote the inner product in $\B^2$ by $(\cdot,\,\cdot)_2$ and by  
 $\E_1,\;\E_2$ we denote
 the projections on the first and the second summands in the orthogonal
 sum $\B^2$, defined by  \eqref{b7}. In the coordinate form we will write for any 
 $h\in \B^2$:
 $
 h = \left\{\E_1 h,\,\E_2h\right\}.
 $
A nonlinear generalization of the theory of self-adjoint extensions 
of linear symmetric operators  
was carried out in the 
papers of the author \cite{Zl-1}-\cite{Zl-5}. 
Before to expound this generalization, we introduce a definition 
\cite{Zl-2}: 
\begin{definition}
Let $F$ be an operator acting in a real Hilbert space $\A$. We 
call it {\it graphically smooth}, if its graph ${\rm graph}(F)$ 
is such a $C^1$-submanifold of $\B^2$, that 
for any $x\in\Rd(F)$ the tangent space 
$T_\zeta({\rm graph}(F))\;(\zeta=\{x,F(x)\})$ is the graph 
of some operator $F^\prime(x)$. We call this operator 
{\it graphic derivative} of $F$ at the point $x$. If 
for any $x\in\Rd(F)$ the operator $F^\prime(x)$ 
is symmetric (self-adjoint), then we call the operator 
$F$ {\it locally symmetric (locally self-adjoint)}.   
\end{definition}
\begin{remark}\label{reb1}
	We mean a smooth submanifold $X$ of some smooth 
	manifold $Y$ as a range of an injective immersion of a smooth 
	manifold $\widehat X$ into $Y$ (\cite{Zl-2}, Definition 1.1). It is
	a submanifold in a generalized sence. When the immersion is a 
	homeomorphism on its range $X$ in the topology induced on the latter
	by the $Y$-topology, then we have the submanifold $X$
	in the commonly accepted sense. According to our terminology, in this
	case it is called {\it the regular submanifold}.      
\end{remark}

We considered in \cite{Zl-2}, \cite{Zl-3} 
a nonlinear operator $F_0$ defined in a real Hilbert space $\A$, 
whose domain is a translated linear set: 
$Dom(F_0)=D_0+\{x_0\}$, where $D_0$ is a dense linear set in $\A$, 
$x_0$ is a distinguished point in $\A$. In applications the space 
$L_2(\Omega)$ plays usually the role of the space $\A$, where 
$\Omega$ is a domain in the space $\R^m$; as the linear set $D_0$ 
participates often the set $C_0^\infty(\Omega)$ of smooth functions 
with compact supports in $\Omega$; $F_0$ is a nonlinear differential 
operator. We suppose that at each point $x\in\Rd(F_0)$ the operator 
$F_0$ has the Gateaux derivative $F^\prime_0(x)$ along $\Rd(F_0)$ 
and the latter is a linear symmetric operator in $\A$ with 
$\Rd(F_0)=D_0$. We considered the following operator 
fields: 
\begin{equation*}
\{\bar F^\prime_0(x)\}_{x\in D_0},\;\;
\{F^\star_0(x)\}_{x\in D_0},
\end{equation*}    
where $\bar F^\prime_0(x)$ is the closure of the operator 
$F^\prime_0(x)$. It turns out, that under some conditions these 
fields are integrable in some sense \cite{Zl-2}. Their antiderivatives 
$\hat F,\;\tilde F$, which satisfy the conditions 
$\hat F(x_0)=\tilde F(x_0)=F_0(x_0)$ are extensions of the operator 
$F_0$ and $F_0\subseteq\hat F\subseteq\tilde F$. We have called 
the operators $\hat F,\;\tilde F$ the {\it minimal} and 
{\it maximal} extensions of the operator $F_0$ respectively. 
They may be considered as nonlinear generalizations of the closure 
$\bar F_0$ for a linear symmetric operator $F_0$ and of the adjoin 
operator $F_0^\star$ to the latter. If $F_0$ is a nonlinear 
differential operator (ordinary or partial elliptic), then its 
maximal extension $\tilde F$ is the maximal differential operator 
with the generalized derivatives (by Sobolev) \cite{Br}, defined 
by the given differential operation in the space $L_2(\Omega)$. 
Notice that under some conditions the operator 
$\tilde F$ is graphically smooth and the operator $\hat F$ is 
locally symmetric.

In the papers \cite{Zl-4}, \cite{Zl-5} the following problem was posed and 
solved: to describe the all locally self-adjoint extensions $F$ of 
the operator $\hat F$, which are restrictions of the 
operator $\tilde F:\;\hat F\subseteq F\subseteq\tilde F$. 
This description was carried out in terms of symplectic differential 
geometry (\cite{Tr}, Chapt. VII). Let us illustrate this for a simple linear ordinary 
differential operator $L_0$ in the space $\A=L_2[0,1]$, defined 
on the set  
\begin{equation*}
\Rd(L_0)=C_0^\infty(0,1)=\{u\in C^\infty(0,1)|\; {\rm supp}
(u)\subset(0,1)\}
\end{equation*}    
by the operation $lu=-\frac{d^2u}{dt^2}$. The operator $L_0$ is 
symmetric. The problem of a description of the all self-adjoint 
boundary value problems for the differential equation 
$l(u)=f$ is equivalent to a description of the all self-adjoint 
extensions of the operator $L_0$. Let $\bar L_0$ be the closure of 
the operator $L_0$. It is clear that $\bar L_0\subseteq L_0^\star$ 
and that every self-adjoint extension $L$ of the operator $L_0$ 
has the property: $\bar L_0\subseteq L\subseteq L_0^\star$. It is 
known \cite{Nm} that the operator $L_0^\star$ is the maximal differential 
operator defined in the space $\A$ by the operation $l$, i.e.
\begin{equation}\label{a1}
\Rd(L_0^\star)=\{u\in\A|\;u,u^\prime\in{\rm Abc}(0,1), lu\in\A\}, 
\end{equation}
where ${\rm Abc}(0,1)$ is the set of the all absolutely continuous 
on $(0,1)$ functions. Furthermore, $L_0^\star u=lu$ for any 
$u\in\Rd(L_0^\star)$. On the other hand, $\Rd(\bar L_0)$    
is described in the following manner:
\begin{equation}\label{a2}
\Rd(\bar L_0)=\{u\in\Rd(L_0^\star)|\;u_0=u_0^\prime=
u_1=u_1^\prime=0\},
\end{equation}
where
\begin{equation}\label{a3}
u_0=u(0+),\;
u_0^\prime=u^\prime(0+),\;
u_1=u(1-),\;
u_1^\prime=u^\prime(1-)
\end{equation}
(these limits exist for any $u\in\Rd(L_0^\star)$ \cite{Nm}).
 Let us light up the above problem of description of 
self-adjoint extensions with the point of view of the symplectic 
geometry.  Consider on $\A^2$ the bilinear form: 
$j(u,v)=(\E_1u,\E_2v)-(\E_2u,\E_1v)$. It can be considered as a closed differential 2-form on $\B^2$, since $j=d\alpha$, where $\alpha$ is the following 1-form on $\B^2$: 
\begin{equation*}
\alpha(\zeta)u:=-(\E_2\zeta,\,\E_1 u)\quad (\zeta,\,u\in\B^2).
\end{equation*}
Since the form $j$ is non-degenerate, it is a symplectic 
form   (\cite{Tr}, Chapt. VII) on the space $\A^2$. It is known that the graphs 
$\Gamma_0={\rm graph}(\bar L_0),\;
\Gamma^\star={\rm graph}(L_0^\star)$ are respectively 
an isotropic and 
a coisotropic subspace of the symplectic space $(\A^2,j)$ and 
a linear operator $L$ acting in $\A$ is self-adjoint 
if and if ${\rm graph}(L)$ is a Lagrangian subspace of 
$(\A^2,j)$ \cite{Zl-4}. Therefore our problem is equivalent to the 
following one: to describe all the Lagrangian subspaces $\Gamma$ of 
$(\A^2,j)$, which contain the isotropic subspace $\Gamma_0$ and 
simultaneously are contained in the coisotropic subspace 
$\Gamma^\star$. These 
subspaces are parameterized by the Lagrangian subspaces 
$\tilde\Gamma$ of the symplectic space $(\C,\tilde j)$, 
where $\C=\Gamma^\star\ominus\Gamma_0$ and $\tilde j$ is the 
lifting of the form $j$ from $\A^2$ into $\C$ by means of the 
natural embedding: $i_\C:\;\C\rightarrow\A^2$. In other words, the 
operator $L$ is a desired self-adjoint extension of $L_0$ if and only 
if ${\rm graph}(L)=\Gamma_0\oplus\tilde\Gamma$ for some Lagrangian 
subspace $\tilde\Gamma$ of $(\C,\tilde j)$. Let us show   
that $\tilde\Gamma$ plays a role of a boundary condition. Using 
the descriptions \eqref{a1}, \eqref{a2}, \eqref{a3} of the domains 
$\Rd(L_0^\star),\;\Rd(\bar L_0)$ and the Green formula:
\begin{equation*}
(lu,v)-(u,lv)=(u^\prime v-uv^\prime)_1-(u^\prime v-uv^\prime)_0\;\;
\forall u,v\in\Rd(L_0^\star) 
\end{equation*}
we can show that the symplectic space $(\C,\tilde j)$ is 
symplectically isometric to the symplectic space  
$(\hat\C,\hat j)$, where 
\begin{equation*}
\hat\C=\{\eta_u\in\R^4:\;\eta_u=(u_0,u_0^\prime,u_1,u_1^\prime),\;
u\in\Rd(L_0^\star)\}
\end{equation*}
and $\hat j(\eta_u,\eta_v)=
(u^\prime v-uv^\prime)_1-(u^\prime v-uv^\prime)_0$. Then the choice 
of a Lagrangian subspace $\hat\Gamma$ of   
$(\hat\C,\hat j)$ is equivalent to the choice of self-adjoint 
boundary conditions.

In the present paper we carry out 
further development of the theory worked out in the 
papers \cite{Zl-4}, \cite{Zl-5} with the purpose to obtain a description 
of locally self-adjoint extensions of the above mentioned nonlinear 
locally symmetric operators in terms of some abstract boundary 
conditions. 

We shall use terminology and notations from the book \cite{Lng} for linear, multilinear and smooth structures and vector bundles.

\begin{remark}
The following question appears: why we restrict ourselves to the case of a real Hilbert space $\B$ ? It turns out that in  a complex Hilbert space the class of graphically smooth locally self-adjoint operators is reduced to the rather small class of ones,  whose graphs are piecewise-affine submanifolds in $\B^2$. This fact is based on the following claim, which follows from Theorem 2.0 of \cite{Gl-Vid}: an analytic function $\Phi(z)$ defined in a domain $G\subseteq\mathbf{C}$, whose values are unitary operators, must be constant. Indeed, the graph of a graphically smooth operator $F$ in a complex Hilbert space $\B$ should be an analytic submanifold of $\B^2$. For any $\zeta\in\mathrm{graph}(F)$ consider the Cayley transform of the derivative $F^\prime(\E_1\zeta)$:
\begin{equation*}
V(\zeta)=\big(F^\prime(\E_1\zeta)-iI\big)\big(F^\prime(\E_1\zeta)+iI\big)^{-1},
\end{equation*}
which is a unitary operator (\cite{Nm}, Chapt. IV, Sect. 14, $n^o$ 3). It is possible to show that the function $V:\,\mathrm{graph}(F)\rightarrow L(\B)$ is analytic in the uniform operator topology. Then using the claim, formulated above, it is possible to show that for any point $\zeta_0\in\mathrm{graph}(F)$ there is its neighborhood $U(\zeta_0)\subset\mathrm{graph}(F)$ such that $F^\prime(\E_1\zeta)\equiv const$ in $U(\zeta_0)$.
\end{remark}
 
\section{\bf Preliminaries}
\setcounter{equation}{0}

\subsection{Maximal and minimal extensions}

 As in \cite{Zl-2}, we consider an operator $F_0$, acting in a real
Hilbert space $\B$, whose domain of definition $\Rd(F_0)$ is
a translated lineal:
\begin{equation*}
\Rd(F_0) = D_0 + \lbrace x_0\rbrace,
\end{equation*}
where $D_0$ is a dense lineal in $\B$, $x_0$ is a distinguished point
in $\B$. We will suppose that the operator $F_0$ satisfies the conditions 
A)-F) of \cite{Zl-2}. Briefly speaking it means, that the operator $F_0$
has at each point $x\in\Rd(F_0)$ the Gateaux derivative 
$F_0^\prime(x)$ along $\Rd(F_0)$ and the latter is a linear 
symmetric operator in $\B$ with $\Rd(F_0^\prime(x)) = D_0$. Denote briefly by $\D(x)$ the closure $\overline{F_0^\prime(x)}$ of the operator $F_0^\prime(x)$.
Furthermore, we imposed the following conditions on the operator field 
$\lbrace\D(x)\rbrace_{x\in D_0}$:

(a) this field admits an extension on the whole $\B$; for the operators 
of the extended field the previous notation is preserved:
$(\D(x))_{x\in \B}$;

(b) there exists a norm $\|.\|_+$, defined on $D_0$, such that the 
completion $\B_+$ of the lineal $D_0$ with respect of this norm is 
continuously imbedded in $\B$ and for any $x\in \B$ the norm 
\begin{equation*}
\|v\|_x = (\|\D(x)v\|^2 + \|v\|^2)^{1\over 2} 
\end{equation*}
is equivalent to the norm $\|v\|_+$; the last fact implies the 
relation:
\begin{equation*}
\forall\,x\in \B\quad
\Rd(\D(x)) = \B_+;
\end{equation*}

(c) the following relation holds:
\begin{equation*}
\D(\cdot)\in C(\B,\,L(\B_+;\B)).
\end{equation*}

We have constructed in \cite{Zl-2} a so called maximal extension of 
the operator $F_0$, which is the conjugate operator $F_0^*$ 
to $F_0$ in the linear case. Let us remind briefly such construction.
For each fixed $x\in\Rd(F_0)$ we consider the following operator
$\Phi^\prime(x)$, which is a member of the class $L(\B;\B_+^*)$:
\begin{equation}\label{b1}
\forall\,h\in \B,\,v\in \B_+\,:\quad
\Phi^\prime(x)h(v) = (\D(x)v,\,h).
\end{equation}
The condition (c) imply, that the operator field 
$\Phi^\prime(x)_{x\in \B}$ is continuous with respect to the operator 
norm. It turns out that this field is integrable (\cite{Zl-2}, Lemma 
2,1).
Let $\Phi$ be its antiderivative, defined by the relation:
\begin{equation}\label{b2}
\Phi(x) = f + \int_{x_0}^x \Phi^\prime(u)du,
\end{equation}
where $f$ is a functional of the form:
\begin{equation}\label{b3}
f(v) = (v,\,F_0(x_0)),
\end{equation}
which belongs to $\B^*\subset \B_+^*$. It is clear that:
\begin{equation}\label{b4}
\Phi\in C^1(\B,\B_+^*).
\end{equation}
We consider the set $\widetilde D = \Phi^{-1}(\B^*\cap\I(\Phi))$,
i.e.
\begin{equation}\label{b5}
\widetilde D = \lbrace x\in \B: \exists\,y\in \B,\,\forall\,v\in \B_+\,:\quad
(\Phi(x))(v) = (v,\,y)\rbrace.
\end{equation}
Since $\B_+$ is dense in $\B$, then in definition \eqref{b5}  
an unique $x\in\widetilde D$ corresponds to each $y\in \B$. Thus we have
an operator $\widetilde F$ acting in $\B$ with 
$\Rd(\widetilde F) = \widetilde D$, defined by the condition:
\begin{equation}\label{b6}
\forall\,v\in \B_+\quad
(\Phi(x))(v) = (v,\,\widetilde F(x)).
\end{equation}
It turns out that the operator $\widetilde F$ is an extension of the 
operator $F_0$ (\cite{Zl-2}, Lemma 2.2). We have called this extension
{\it the maximal extension} of the operator $F_0$.  
According to the definition \eqref{b6} of the operator $\widetilde F$, 
its graph
\begin{equation}\label{b8}
\M =\T(\widetilde F)
\end{equation}
is the set of zeros of an operator $\Theta$, which is a member of the 
class $C^1(\B^2, \B_+^*)$ and is defined by the relation:
\begin{equation}\label{b9}
(\Theta(\zeta))(v) = (\Phi(x))(v) - (v,\,y),
\end{equation}
where $v\in \B_+\quad \zeta = {x, y}\in \B^2$. By Theorem 2.1 of \cite{Zl-2}
the set $\M$ is a regular $C^1$-submanifold of $\B^2$. So the
operator $\widetilde F$ is closed and graphically smooth (\cite{Zl-2},
Definition 1.2) and at each point $x\in\widetilde D$ its graphic 
derivative $\widetilde F^\prime(x)$ is of the form:
\begin{equation}\label{b8}    
\widetilde F^\prime(x) = (\D(x))^*
\end{equation}
(\cite{Zl-1}, Theorem 2.1).

By Theorem 4.1 of \cite{Zl-2} the domain of definition $\widetilde D$ 
\eqref{b5} of the operator $\widetilde F$ is $\B_+$-saturated, i.e.
\begin{equation}\label{b11}
\widetilde D = \widetilde D + \B_+.
\end{equation}
Then we can define the following extension of the operator $F_0$:
\begin{equation}\label{b12}
\widehat F = \widetilde F\vert_{\left\{x_0\right\}+\B_+}.
\end{equation}
We have called it {\it the minimal extension} of the operator $F_0$. 
It may be considered as a nonlinear generalization of the closure
of the operator $F_0$. Indeed, by Lemma 4.1 of \cite{Zl-2} the operator
$\widehat F$ is graphically smooth and at each point 
$x\in\Rd(\widehat F$ its graphic derivative 
$\widehat F^\prime(x)$ coincides with $\D(x)$. Let us remind
that the last operator is the closure of $F_0^\prime(x)$.
Since for each $x\in H$ the operator $\D(x)$ is symmetric,
then the operator $\widehat F$ is {\it locally symmetric} 
(\cite{Zl-3}, Definition 6.1).

In \cite{Zl-2} we have defined on the manifold $\M$ \eqref{b8} the following
Abelian group of mappings:
\begin{equation}\label{b13}
\G = \lbrace G(\cdot,v)\rbrace_{v\in \B_+},
\end{equation}
where
\begin{equation}\label{b14}
\forall\,\zeta\in\M,\quad \forall\,v\in \B_+ \quad
G(\zeta,v) = \lbrace\E_1\zeta + v,\,
\widetilde F(\E_1\zeta + v)\rbrace.
\end{equation}
In view of the property \eqref{b11}, the right part of \eqref{b14} is well 
defined. This group can be considered as a representation
of the additive group $\B_+$, i.e.
\begin{equation}\label{b15}
\forall\,\zeta\in\M,\quad \forall\,v_1,v_2\in H_+ \quad
G(G(\zeta,v_1),v_2) = G(\zeta, v_1+v_2).
\end{equation}

It turns out, that the manifold $\M$ is a regular 
$C^2-$ submanifold of $\B^2$ and
\begin{equation}\label{b16}
G(\cdot,\cdot)\in C^1(\M\times \B_+,\,\M),
\end{equation}
if the operator field $\lbrace\D(x)\rbrace_{x\in \B}$ 
satisfies the 
following additional conditions (\cite{Zl-2}, Lemma 4.2, \cite{Zl-3},
Theorem 5.1):

(d) if the mapping $\D(\E_1\zeta)$ is considered
as acting from $\M$ into \\$L(\B_+;\B)$, then at each 
$\zeta\in \M$ it has the linear Gateaux derivative\\
$(\D(\E_1\zeta)\cdot)^\prime\cdot$;

(e) for any fixed $\zeta\in \M,\;\;v\in \B_+$ the operator
$(\D(\E_1(\zeta)v)^\prime\cdot$, which is defined on
$\E_1(T_\zeta(\M))$, admits an extension on all
of $\B$, which we denote by 
$\D^\prime(\E_1(\zeta)(v,\cdot)$; 

(f) $\D^\prime(\E_1\zeta)(\cdot,\cdot)$ is a locally uniformly
continuous mapping from $\M$ into $L(\B_+,\B;\B)$.

\subsection{Locally self-adjoint extension and the foliation generated by \\the group $\G$}

 We have called a graphically smooth 
operator $F$  {\it locally self-adjoint}, if its graphic derivative $F^\prime(x)$ at any point
$x\in \Rd(f)$ is a self-adjoint operator (\cite{Zl-3},
Definition 6.1). In the papers \cite{Zl-4}, \cite{Zl-5} the following
problem was raised and solved: to describe all locally self-adjoint
extensions $F$ of the operator $F_0$, which satisfy the condition:
\begin{equation}\label{b17}
\widehat F_0\subset F\subset \widetilde F.
\end{equation}
The description mentioned above was carried out in terms of the
{\it symplectic differential geometry} (\cite{Tr}, Chapt VII). We consider the space
$\B^2$ \eqref{b7} as a symplectic manifold $(\B^2,\,j)$ with the symplectic
form $j$, desribed in Introduction, i.e.,
\begin{equation*}
\forall\, \zeta\in \B^2,\quad \forall\,h_1, h_2\in T_\zeta(\B^2)
(\equiv \B^2)\quad
j(\zeta)(h_1,h_2) = (Jh_1,\,h_2)_2,
\end{equation*}
where the operator $J$ is defined by the relation: 
\begin{equation*}
\forall\,h\in \B^2\quad
Jh = \lbrace-\E_2 h_2,\,\E_1 h\rbrace.
\end{equation*}
This operator is sqew-self-adjoint and unitary:
\begin{equation*}
J^* = -J,\quad J^2 = I.
\end{equation*} 
Notice that we identify the tangent bundle $T(\B^2)$ with 
$\B^2\times \B^2$.

A graphically smooth operator $F$ is {\it locally symmetric (locally
self-\\-adjoint)}, iff $\T(F)$ is an {\it isotropic (Lagrangian)}
submanifold of the synplectic manifold $(\B^2,\,j)$ (\cite{Zl-4},
Proposition 1.1).

The equality \eqref{b8} implies, that $\M$ is a {\it coisotropic}
submanifold of \\$(\B^2,\,j)$ \cite{Zl-4}.

The above raised problem on locally self-adjoint extensions my be 
reformulated in the following manner: to describe all Lagrangian
submanifolds $\K$ of $(\B^2,\,j)$, such that
\begin{equation}\label{b18}
\Gamma_0\subset\K\subset\M,
\end{equation}
where
\begin{equation}\label{b19}
\Gamma_0 = \T(\widehat F).
\end{equation}
Such manifolds $\K$ are the graphs of desired locally self-adjoint
extensions $F$ of the operator $F_0$, satisfying the condition
\eqref{b17}. This problem was solved in \cite{Zl-5} also in a more restricted
sence. We have described so called {\it regularly} locally self-adjoint
extensions, i.e. such ones, whose graphs are regular submanifolds
of $\B^2$ 
(\cite{Zl-3}, Remark 2.1).

Let us remind the main statements of \cite{Zl-4} concerning above 
mentioned extensions. In their construction an important role play
the orbits $\Gamma_\zeta$ of the group $\G$ \eqref{b13}, \eqref{b14}), 
which have the form:
\begin{equation*}
\forall\,\zeta\in\M\quad
\Gamma_\zeta = \T(\widehat F_\zeta),
\end{equation*}
\begin{equation}\label{b20}
\widehat F_\zeta = \widetilde F\vert_{\lbrace\E_1\zeta\rbrace 
+ \B_+}.
\end{equation}
In particular
\begin{equation}\label{b21}
\widehat F = \widehat F_{\zeta_0},
\end{equation}
where
\begin{equation*}%\label{b22}
\zeta_0 = \lbrace x_0,\,F_0(x_0)
\end{equation*}
(see (2.12)). Furthermore:
\begin{equation*}
\Gamma_0 = \Gamma_{\zeta_0}
\end{equation*}
(see \eqref{b19}).

The orbits $\Gamma_\zeta$ have the property:
\begin{equation*}
\forall\,\zeta_1,\zeta_2\in\M\quad
{\rm either}\quad \Gamma_{\zeta_1}\cap\Gamma_{\zeta_2} = \emptyset
\quad{\rm or}\quad \Gamma_{\zeta_1} = \Gamma_{\zeta_2}.
\end{equation*}
This means that the family of these orbits 
\begin{equation}\label{b24}
\F = \lbrace\Gamma_\zeta\rbrace_{\zeta\in\M}.
\end{equation}
forms a partition of the manifold $\M$. It turns out that under the
conditions (a)-(f) this partition forms a $C^1$-{\it foliation} on
$\M$ in the sence of R.S. Palais, i.e. it is a local 
$C^1$-bundle (\cite{Pal}; \cite{Zl-4}, Definition 2.1). A {\it fibered
chart}, corresponding to above foliation near any point
$\zeta_*\in\M$, may be chosen in the form:
\begin{equation}\label{b25}
(U,\psi^{-1}, \B_+\times B),
\end{equation}
where
\begin{equation}\label{b26}
U = \psi(V\times W),
\end{equation}
$V$ is a neighborhood of 0 in $\B_+$,
\begin{equation}\label{b27}
W = \phi(\widehat W)
\end{equation}
and
\begin{equation}\label{b28}
(\widehat W, \phi, B),\quad \zeta _*\in W
\end{equation}
is a chart at the point $\zeta _*$ for a regular $C^1$-submanifold 
$\C$ of $\M$ with a model Banach space $B$. This submanifold
complements each leaf $\Gamma_\zeta$ in 
$\M$. This  means that if $\C\cap\Gamma_\zeta\neq\emptyset$, then
\begin{equation*}
\forall\,\xi\in \C\cap\Gamma_\zeta\quad
T_\xi(\Gamma_\zeta)\cap T_\xi(\C) = \{0\}\quad\mathrm{and}\quad T_\xi(\Gamma_\zeta)+ T_\xi(\C)= T_\xi(\M).
\end{equation*}
The mapping $\psi$, which participates
in the chart \eqref{b25} has the form:
\begin{equation*}
\forall\,\eta\in W,\,v\in V\quad
\psi(v,\eta) = G(\phi^{-1}(\eta),v).
\end{equation*}
Hence the leafs $\Gamma_\zeta$ are $C^1$-submanifolds of $\M$
with the model space $H_+$. We called linear connected components
of sets $\psi(V\times W)\cup\Gamma_\zeta$ {\it sections} of the 
fibered chart \eqref{b25}. In our case such sections have the form:
\begin{equation*}
\sigma_\zeta = \psi(V\times {\eta}),\; \eta\in W.
\end{equation*}

We dealt with so called {\it regular foliation} (\cite{Pal}; 
\cite{Zl-4}, Definition 2.3). This means that near each point 
$\zeta\in \M$ there exists a so called {\it regular
fibered chart} of the form \eqref{b25}, such that each leaf $\Gamma_\zeta$
of the foliation intersects $U$ in at most one section
$\sigma_\zeta$ .
The following conditions ensure the regularity of the foliation
$\F$, defined by \eqref{b24} (\cite{Zl-4}, Theorem 2.1):

(g) the coercive estimate holds for the maximal extension 
$\widetilde F$ of the operator $F_0$:
\begin{equation*}
\forall\,x\in\Rd(\widetilde F),\,v\in \B_+\quad
\|\widetilde F(x+v) - \widetilde F(x)\|^2 + \|v\|^2
\geq \gamma(\|x\|, \|v\|_+),
\end{equation*}
where $\gamma(y,z)$ is a continuous function in the domain 
$y\geq 0,\;z\geq 0$ and satisfies the conditions:

1) $\gamma(y,z)>0$, when $z\neq 0$;

2) $\gamma(y,0)=0$;

3) for any fixed $y\geq 0$ the function $\gamma(y,z)$ increases with
respect to $z$.

If the foliation $\F$ is regular, then it may be endowed with the 
structure of a $C^1$-manifold, whose topology coincides with the
quotient topology, generated by the partition of $\M$ into 
leafs of $\F$. This structure is defined in the following manner.
Consider near any point $\zeta\in\M$ a regular fibered 
chart of the form \eqref{b25}-\eqref{b26}. Let $\Pi_\F$ be
the natural projection of $\M$ on the partition $\F$, i.e            
\begin{equation}\label{b31}  
\Pi_\F(\zeta) = \Gamma_\zeta.
\end{equation}
We set
\begin{equation}\label{b32}
\widehat U = \Pi_{\F}(U).
\end{equation}
Then the submanifold
\begin{equation}\label{b33}
\widehat W = \psi({0}\times W)
\end{equation}
intersects each leaf $\Gamma\in\widehat U$ at an unique point
\begin{equation}\label{b34}
\xi = \pi(\Gamma).
\end{equation}
Consider the mapping
\begin{equation}\label{b35}
\widetilde\phi = \phi\cdot\pi,
\end{equation}
where $\phi$ is the mapping, participating in the chart \eqref{b28}. 
It turns out that the charts of the form
\begin{equation}\label{b36}
(\widehat U, \widehat \phi, B)
\end{equation}
(see \eqref{b31}, \eqref{b32}, \eqref{b34}, \eqref{b35}) 
constitute a $C^1$-atlas
on $\F$ \cite{Zl-4}.

Since in view of \eqref{b8}, \eqref{b20},
\begin{equation}\label{b37}
\forall\,\zeta,\,\xi\in \M\quad
(\widehat F_\zeta)^\prime(\E_1\xi) = 
\D(\E_1\xi),
\end{equation}
then the leafs $\Gamma_\zeta$ are isotropic submanifolds of the 
symplectic manifold \\$(\B^2,\,j)$. Furthermore, the pair
\begin{equation}\label{b38}
(\widehat W,\,\omega),
\end{equation}
where
\begin{equation}\label{b39}
\omega = i_{\widehat W}^*(j),
\end{equation}
is a symplectic manifold, since at each $\zeta\in \widehat W$ the 
leaf $\Gamma_\zeta$ complements the submanifold $\widehat W$, defined by \eqref{b33},           
in the coisotropic submanifold $\M$. Hence the pair
\begin{equation}\label{b40}
(\widehat U,\,\widehat\omega),
\end{equation}
where
\begin{equation}\label{b41}
\widehat\omega = \pi^*(\omega)
\end{equation}
is a symplectic manifold (see \eqref{b32}, \eqref{b34}). Thus, the foliation
$\F$ may be considered as a symplectic manifold in the local
sence, described above.

Let us formulate the main result of the paper \cite{Zl-4} (Theorem 3.1):
\begin{proposition}\label{prb1}
Assume the foliation $\F$ is
regular and the conditions (a)-(f) are satisfied. 
A graphically smooth operator $F$ is regularly locally self-adjoint
and and satisfies the condition \eqref{b17}, if the following
conditions are satisfied:

1) $\T(F)$ is a $\F$-saturated set, i.e. for any 
$\zeta\in \T(F)$ the lief $\Gamma_\zeta\subset\T(F)$;

2) for any leaf $\Gamma\in\Pi_{\F}(\T(F))$ (see \eqref{b31})
there exists its neighborhood $\widehat U\subset\F$ of 
the form \eqref{b32}, such that the set 
\begin{equation*}
\K = \widehat U\cap\Pi_{\F}(\T(F))
\end{equation*}
is a regular Lagrangian submanifold of the symplectic manifold 
$(\widehat U,\,\widehat\omega)$.

If $\T(F)$ is a closed $C^2$-submanifold in $\M$, then
above conditions are necessary.
\end{proposition}

We will describe a specific class of regularly locally self-adjoint
extensions of the operator $F_0$. For this aim we need the following
definition:
\begin{definition}\label{deb1}

Let $N$ be a subset of $\M$. We call 
the set
\begin{equation*}
\Pi_{\F}^{-1}(\Pi_{\F}(N)) =
\bigcup_{\zeta\in N}\,\Gamma_\zeta
\end{equation*}
the $\F$-{\it saturation of the set} $N$ and denote it by  
$N_s$.

We call $N$ $\F$-{\it saturated}, if $N_s = N$.
\end{definition}

The following statement is valid (\cite{Zl-4}, Theorem 3.4):
\begin{proposition}\label{prb2}
Assume that the foliation $\F$ 
is regular and the conditions (a)-(f) are satisfied. Let
$\widehat W$ be a regular $C^1$-submanifold of $\M$ such that
each leaf $\Gamma$ of the foliation $\F$ intersects  
the latter in at most one point and complements it in $\M$ at each intersection point (if it exists). Then in order
an operator $F$ would be a regularly self-adjoint extension of $F_0$
it is sufficient that
\begin{equation*}
\T(F) = L_s,
\end{equation*}
where $L$ is a regular Lagrangian submanifold of the symplectic 
manifold \\$(\widehat W,\,\omega)$, defined by \eqref{b38}, \eqref{b39}.
\end{proposition} 

\begin{definition}\label{deb2}
We call the operator $F$, described in 
Proposition \ref{prb2}, the locally self-adjoint extension, defined
by the Lagrangian submanifold $L$ of the symplectic 
manifold $(\widehat W,\,\omega)$.
\end{definition}

\section{\bf Abstract boundary conditions}
\setcounter{equation}{0}

The problem of a description of above mentioned locally 
self-adjoint extensions of the operator $F_0$ in terms of so called 
{\it abstract boundary conditions} can be formulated explicitly
in the following manner: what properties the mapping
$\Theta:\; \M\rightarrow E$ ($E$ is a Banach space) must
have in order the set of its zeros
\begin{equation*}\label{c1}
\K = \{\zeta\in\M:\;
\Theta(\zeta) = 0\}
\end{equation*}
would be the graph of the desired extension ?

\subsection{Symplectomorphic property of the group $\G$}

Before to solve above problem, we will establish a symplectomorphic property
of the mappings $G(\cdot,v)$ of the group $\G$.
\begin{proposition}\label{prc1} 
Under the conditions (a)-(f) each 
mapping $G(\cdot,v)\;(v\in \B_+)$ realizes a $C^1$-symplectomorphism 
of the coisotropic $C^2$-submanifold of the symlectic manifold
$(\B^2,\,j)$ into itself, i.e.
\begin{eqnarray}\label{c2}
&&\forall\, v\in \B_+\;\forall \zeta\in\M,\;
\forall\,h_1,h_2\in T_\zeta(\M)\quad\nonumber\\
&&(JG_\zeta(\zeta,v)h_1,\,G_\zeta(\zeta,v)h_2)_2 = 
(Jh_1,\,h_2)_2.
\end{eqnarray}
\end{proposition}
\begin{proof}
For any fixed $v\in \B_+$ consider the vector field on 
$\B^2$:
\begin{equation}\label{c3}
\xi_v(\zeta) = \{v,\,\D(\E_1\zeta)v\}.
\end{equation}
Definition \eqref{b14} of the mappings $G(\cdot,v)$ and the relation
\begin{equation*}
\forall\, \zeta\in\M\quad
\D(\E_1\zeta)\subset 
\widetilde F^\prime (\E_1\zeta)
\end{equation*}
imply, that the restriction
\begin{equation}\label{c4}
\widehat \xi_v = \xi_v\vert_{\M}
\end{equation}
forms a vector field on $\M$ and the restriction of the flow
$U(t,\zeta,v)$ of the field $\xi_v$ \eqref{c3} on $\M$ coincides
with the flow $\widehat U(t,\zeta,v)$ of the field $\widehat\xi_v$,
which has the property:  
\begin{equation}\label{c5}
\forall\,\zeta\in \M\quad
\widehat U(t,\zeta,v) = U(t,\zeta,v) = G(\zeta, tv).
\end{equation}
In view of the conditions d), e), f), the field $\widehat \xi_v$
belongs to the class $C^1$.

For the proof of our proposition it is sufficient to show, that the
vector field $\xi_v$ is Hamiltonian, i.e. there exists a scalar function
$f\in C^2(\B^2,R)$, such that
\begin{equation}\label{c6}
\forall\,\zeta\in \B^2,\,h\in \B^2\quad
(J\xi_v(\zeta),\,h) = df(\zeta)h.
\end{equation}
Indeed, assuming that the relation \eqref{c6} holds, consider the Lie
derivative $L_{\widehat \xi_v}(\widehat j)$ of the 2-form $\widehat j = i_{\M}^*j$
along the field $\widehat \xi_v$ defined by \eqref{c4} (\cite{Tr}, Chapt. VII, Sect 2). Putting 
$\widehat f = f\vert_{\M}$, we obtain according to the 
definition of the Lie derivative:
\begin{eqnarray*}
&&\forall\,\zeta\in \M\quad
{d\over dt}(\widehat U(t,\zeta,v)^*\widehat j)\vert_{t=0} =
L_{\widehat \xi_v}(\widehat j) =\nonumber\\
&&d(\widehat j\vee \widehat \xi_v) +D(\widehat j)
\vee \widehat \xi_v = d(\widehat j\vee \widehat \xi_v) =
d((J\widehat \xi_v,\,\cdot)_2) = d(d\widehat f(\zeta)\cdot) = 0,
 \end{eqnarray*}
hence the flow $\widehat U(t,.,v)$ is a symplectomorphism for each fixed $t\in\R$. Then 
desired property \eqref{c2} follows from \eqref{c5}.

We now turn to the proof of the relation \eqref{c6}. The latter is 
equivalent to the fact that the vector field 
\begin{equation}\label{c7}
\eta_v(\zeta) = J\xi_v(\zeta)
\end{equation}
is integrable, i.e. for any closed contour $C\subset H^2$ 
\begin{equation}\label{c8}
\int_C(\eta_v(\zeta),\,d\zeta)_2 = 0.
\end{equation}
Let $\Lambda$ be the following lineal in $\B^2$:            
\begin{equation}\label{c9}
\Lambda =\Rd(F_0)\oplus \B,
\end{equation}
which is dense in $\B^2$:
\begin{equation}\label{c10}
{\rm cl}(\Lambda) = \B^2.
\end{equation}
Suppose $\zeta_0,\zeta_1\in \Lambda$ and $L_{0,1}$ is 
the segment of the straight, connecting the points $\zeta_0,\;\zeta_1$:
\begin{equation*}%\label{c11}
L_{0.1} = (\zeta\in \B^2: \zeta = \zeta(t) = 
(1-t)\zeta_0 + t\zeta_1,\; 0\leq t\leq 1).
\end{equation*}
Consider the integral:
\begin{equation*}%\label{c12}
I_{0,1} = \int_{L_{0,1}}(\eta_v(\zeta),\,d\zeta)_2  
\end{equation*}
(see \eqref{c3}, \eqref{c7}). We have:
\begin{eqnarray*}%\label{c13}
&&\hskip-15mm I_{0,1} = \int_0^1(\D(\E_1(\zeta(t))v,\,
\E_1(\zeta_1-\zeta_0))dt - 
(v, \E_2(\zeta_1-\zeta_0)) =\nonumber\\
&&\hskip-15mm(v,\,\int_0^1 F_0^\prime (\E_1(\zeta_1 -\zeta_0))dt) 
- (v, \E_2(\zeta_1-\zeta_0)) =
\Psi(\zeta_1,v) - \Psi(\zeta_0,v),
\end{eqnarray*}
where
\begin{equation*}%\label{c14}
\Psi(\zeta_,v) = (v,\,F_0(\E_1\zeta) - \E_2\zeta).
\end{equation*}
The last representation 
of the integral $I_{0,1}$ implies:
\begin{equation*}
\int_\Delta(\eta_v(\zeta),\,d\zeta)_2 = 0
\end{equation*}
for any triangular contour $\Delta$ with vertices in $\Lambda$  
\eqref{c9}. From \eqref{c10} and 
the continuity of the field $\eta_v$ we
conclude, that the same equality holds for any triangular contour
$\Delta\subset \B^2$. Then by the Gavurin theorem \cite{Vn}
the relation \eqref{c8} is valid, i.e. the field $\eta_v$ is integrable.
\end{proof}

\subsection{Defect of the operator $F_0$}

Consider a linear closed symmetric operator 
$A_0$ acting in a real Hilbert space $\B$ with a dense domain $\Rd((A_{0})$, the complexification $\B_c = \B\otimes_\R\mathbf{C}$ of 
$\B$ and the natural extension  $A_{0,c}$ of $A_0$ to $\B_c$, defined by $A_{0,c}(x\otimes 
z)=(A_{0}x)\otimes z$ for any $x\in\Rd((A_{0})$ and $z\in\mathbf{C}$. It is known that $A_{0,c}$ is a symmetric operator acting in $\B_c$ with the dense domain $\Rd((A_{0,c})=\Rd((A_{0})\otimes\mathbf{C}$.
We shall consider deficiency index of $A_{0,c}$ (\cite{Nm}, Chapt IV, Sect.14, $n^o$ 7).
The following claim is valid:

\begin{proposition}\label{prdefect}
(i) Defect numbers of the operator $A_{0,c}$ are equal each to other, i.e., its deficiency index has the form $(m,m)$;

(ii) The equality
\begin{equation}\label{c15}
\dim(\T(A_0^*)\ominus \T(A_0)) = 2m.
\end{equation}
is valid.	
\end{proposition}
\begin{proof}
(i) Let $(n_i, n_{-i})$ be the deficiency index of $A_{0,c}$. This means that 
$n_i=\dim_c(\mathcal{N}_i)$ and $n_{-i}=\dim_c(\mathcal{N}_{-i})$, where
$\mathcal{N}_i=\ker(A_{0,c}^\star-iI)$ and $\mathcal{N}_{-i}=\ker(A_{0,c}^\star+iI)$. Here $\dim_c$ denotes the complex dimension of subspaces in  $\B_c$. It is not difficult to show that $JA_{0,c}^\star=A_{0,c}^\star J$, where $J$ is the conjugation operator in $\B_c$. Hence $\mathcal{N}_i=J\mathcal{N}_{-i}$, i.e., $n_i=n_{-i}$. Claim (i) is proven.

(ii)  By the Neumann formula for $\Rd(A_{0,c}^\star)$ (\cite{Nm}, Chapt IV,
 Sect.14, Theorem 4):
$\Rd(A_{0,c}^\star)=\Rd(A_{0,c}^\star)+\mathcal{N}_i+\mathcal{N}_{-i}$, where the summands are linearly independent. Hence, in view of claim (i), 
\begin{equation*}
\dim_c(\Rd(A_{0,c}^\star))
=2m\,(\mathrm{mod}(\Rd(A_{0,c}))), 
\end{equation*}
therefore $\dim(\Rd(A_{0}^\star))=2m\,(\mathrm{mod}(\Rd(A_{0})))$.
 The last equality implies easily \eqref{c15}. Claim (ii) is proven.
\end{proof}

\begin{definition}
We shall call the number $m$, definef by \eqref{c15}, {\it defect} of the operator $A_0$ and denote it by ${\rm def}(A_0)$. 
\end{definition}
 
We shall deal with a nonlinear operator $F_0$, for which the closure of $\D(x)$ of its graphic
derivative $(F_0^\prime(x)$ has at all points of $\Rd(F_0)$ the same finite 
defect:
\begin{equation}\label{c16}
\forall\,x\in \Rd(F_0)\quad
{\rm def}(\D(x)) = m < \infty.
\end{equation}
It turns out, that the relation \eqref{c16} holds for any  
$x\in \Rd(F_0)$, if it is true for at least one 
$x_0\in \Rd(F_0)$. More explicitly, the following statement
is true:
\begin{proposition}\label{prc2}
If the conditions (a)-(d) are fulfilled and
\begin{equation}\label{c17}
\exists\,x_0\in \Rd(F_0)\quad
{\rm def}(\D(x)) = m < \infty,
\end{equation}
then
\begin{equation}\label{c18}
\forall\,x\in \B\quad
{\rm def}(\D(x)) = m.
\end{equation}
\end{proposition}
\begin{proof}
In view of the condition (c), the subspaces 
$\T(\D(x))$ depend continuously on $x$ with 
respect to the gap metric (\cite{Kat}, Chapt. IV, Theorem 2.14). The subspaces 
$\T((\D(x))^*)$ have the same property, since
\begin{equation*}
\T((\D(x))^*) =
J(\T(\D(x)))^\bot.
\end{equation*}
Therefore the subspaces
\begin{equation*}
\C_x = \T((\D(x))^*)\ominus 
\T(\D(x))
\end{equation*}
depend continuously on $x$ with respect to the gap metric too, 
hence their dimensions are equal to each other. Thus, in view of Proposition \ref{prdefect}, \eqref{c18}
is valid.
\end{proof}

The proved proposition brings us to the following definition:
\begin{definition}\label{dec1}
If for the operator $F_0$ the condition \eqref{c17} 
is fulfilled, then we say, that $F_0$ has a {\it finite defect}
$m$ and write $m = {\rm def}(F_0)$.
\end{definition}

\subsection{Description of abstract boundary conditions in the linear case}

We now turn to description of abstract boundary conditions. Before
we will establish corresponding statement for the linear case.
Let $A_0$ be a closed symmetric operator in $\B$ with a finite 
defect $m$ and $\Theta$ be a linear operator such that     
\begin{equation}\label{c19}
\Theta\in L(\M, \R^m),
\end{equation}
where
\begin{equation}\label{c20}
\M = \T(A_0^*).
\end{equation}
Since $\M$, $\R^m$ are Hilbert spaces, then
\begin{equation}\label{c21}
\Theta ^*\in L(\R^m,\M).
\end{equation}

The following statement describes self-adjoint extensions of 
$A_0$ in terms of abstract boundary conditions:
\begin{lemma}\label{lec1}
The subspace
\begin{equation}\label{c22}
\K ={\rm ker}(\Theta)
\end{equation}
is the graph of a self-adjoint extension $A$ of the operator $A_0$
iff the following conditions are fulfilled:  
\begin{equation}\label{c23}
\I(\Theta) = \R^m;
\end{equation}
\begin{equation}\label{c24}
\I(\Theta^*)\subset \C,
\end{equation}
where
\begin{equation}\label{c25}
\C = \M\cap J(\M),
\end{equation}
and
\begin{equation}\label{c26}
\Theta J\Theta^* = 0.
\end{equation}
\end{lemma}
\begin{proof}
Assume that 
\begin{equation}\label{c27}
\K = \T(A),
\end{equation}
\begin{equation}\label{c28}
A_0\subset A,
\end{equation}
\begin{equation}\label{c29}
A = A^*.
\end{equation}
Consider the subspace
\begin{equation}\label{c30}
\Gamma = \T(A_0).
\end{equation}
Since $A_0$ is a symmetric operator, then it is a restriction of the
operator $A^*$. This means that $\Gamma$ is an isotropic subspace
of the symplectic space $(\B^2,\,j)$, which is contained in the 
coisotropic subspace $\M$, defined by \eqref{c20}:
\begin{equation}\label{c31}
\Gamma\subset\M.
\end{equation}
Furthermore, in view of \eqref{c20}, we have:
\begin{equation}\label{c32}
\M = (J\Gamma)^\bot.
\end{equation}
According to \eqref{c22}, \eqref{c27}, we can rewrite \eqref{c28} 
in the form:
\begin{equation}\label{c33}
\Gamma\subset\ker(\Theta).
\end{equation}
From \eqref{c25}, \eqref{c32} we conclude:
\begin{equation}\label{c34}
\C = \M\cap\Gamma^\bot.
\end{equation}
Then the inclusion \eqref{c33} is equivalent to the following one:
\begin{equation*}
\I(\Theta^*) = \M\ominus\ker(\Theta)
\subset\C.
\end{equation*}
So \eqref{c28} is equivalent to \eqref{c24}-\eqref{c25}. 
From \eqref{c25} we conclude:
\begin{equation}\label{c35}
J(\C) = \C.
\end{equation}

Then the 2-form $\widehat j(u,v) = (\widehat Ju,\,v)_2$, where
\begin{equation}\label{c36}
\widehat J = J\vert_{\C},
\end{equation}
is a symplectic form in the space $\C$. So the pair
$(\C,\;\;\widehat j)$ form a symplectic space. In view 
of \eqref{c32} and the equality \eqref{c34},
\begin{equation}\label{c37}
\dim(\C) = 2m.
\end{equation}
Consider the subspace:
\begin{equation}\label{c38}
L = \K\cap\C.
\end{equation}
Then by \eqref{c34} we have:
\begin{equation}\label{c39}
\K = L\oplus \Gamma.
\end{equation}
Since $\Gamma$ is an isotropic subspace of $(\B^2,\,j)$,
it is the $J$-orthogonal complement to the coisotropic subspace 
$\M$ (see \eqref{c32}) and \eqref{c31}) is valid. Then the representation 
\eqref{c39} implies that the subspace $\K$ is Lagrangian in 
\begin{equation*}
(\B^2,\,j),
\end{equation*}
iff the subspace $L$ \eqref{c38} is Lagrangian in 
$(\C,\;\; \widehat j)$.

We have by \eqref{c22}), \eqref{c38} that
\begin{equation}\label{c40}
L = \ker(\widehat\Theta),
\end{equation}
where
\begin{equation}\label{c41}
\widehat\Theta = \Theta\vert_\C.
\end{equation}
From \eqref{c33}, \eqref{c34} we conclude , that the equality 
\eqref{c23} is 
equivalent to the following one:
\begin{equation}\label{c42}
\I(\widehat\Theta) = \R^m.
\end{equation}
Furthermore, by \eqref{c24}, \eqref{c41} 
$\widehat\Theta^* = \Theta^*$. Then
taking into account \eqref{c35},\,\\\eqref{c36}, 
we obtain that the relation \eqref{c26}
is equivalent to
\begin{equation*}
\widehat\Theta\widehat J\widehat\Theta^* = 0.
\end{equation*}
Since by \eqref{c35} $\widehat J^{-1} = -\widehat J$, the last equality can
be rewritten in the form:
\begin{equation}\label{c43}
\widehat\Theta\widehat J^{-1}\widehat\Theta^* = 0.
\end{equation}
By \eqref{c37} and Proposition 4.1 the relatons 
\eqref{c41}, \eqref{c42}, \eqref{c43} 
are equivalent to the fact, that $L$ is a Lagrangian subspace of the
symplectic space $(\C,\;\;\widehat j)$, i.e. as it was 
mentioned above, $L$ is a Lagrangian subspace of $(\B^2,\,j)$. The
last fact is equivalent to \eqref{c29}. So the relations \eqref{c27}-
\eqref{c29}
are equivalent to \eqref{c23}-\eqref{c26}.
\end{proof}

\subsection{Main results}

We now turn to a description of abstract locally self-adjoint boundary
conditions of our nonlinear operator $F_0$, which we suppose to have
a finite defect $m$ (Definition \ref{dec1}).
\begin{theorem}\label{tec1}
Assume that the conditions (a)-(f) are
fulfilled and the foliation $\F$ is regular. A non-emty subset
$\K\subset\M$ is the graph of a locally 
self-adjoint extension of the operator $F_0$, if the following
conditions are fulfilled:

1) $\K\cap\T(F_0)\neq\emptyset$;

2) there exists a neighborhood $U(\K)$ of $\K$ in $\M$
and a mapping
\begin{equation}\label{c44}
\Theta\in C^2(U(\K),\,\R^m),
\end{equation}
such that $\K$ satisfies the condition 
\begin{equation}\label{c1}
\K = \{\zeta\in\M:\;
	\Theta(\zeta) = 0\}
\end{equation}
and at each point
$\zeta\in U(\K)$ the following conditions are fulfilled:        
\begin{equation}\label{c45}
\I(\Theta^\prime(\zeta)) = \R^m,
\end{equation}
\begin{equation}\label{c46}
\I((\Theta^\prime(\zeta))^*)\subset \C_\zeta,
\end{equation}
where
\begin{equation}\label{c47}
\C_\zeta = T_\zeta(\M)
\cap J(T_\zeta(\M))
\end{equation}
and
\begin{equation}\label{c48}
\Theta^\prime(\zeta) J(\Theta^\prime(\zeta))^* = 0.
\end{equation}
\end{theorem}

\begin{proof}
By conditions \eqref{c44}-\eqref{c45} the set $\K$
is a $C^2$-submanifold of $\M$ and
\begin{equation*}
\forall\,\zeta\in\K\quad
T_\zeta(\K) = \ker(\Theta^\prime(\zeta)).
\end{equation*}
This fact,  conditions \eqref{c45}-\eqref{c48} 
and Lemma \ref{lec1} imply that for any $\zeta\in\K$ the subspace $T_\zeta(\K)$ is the graph of a self-adjoint operator.
This means that $\K$ is the graph of a locally self-adjoint operator
$F$. By the condition 1) there exists a point 
$y_0\in\E_1(\K)\cap\Rd(F_0)$. On the other
hand by the Lemma 3.1 of \cite{Zl-4} the submanifold $\K$ is
$\F$-saturated, hence
\begin{equation*}
\Gamma_{\xi_0} = \T(\widehat F)
\subset \K,\quad
\xi_0 = \lbrace y_0,\,F_0(y_0)\rbrace.
\end{equation*}
Recall that $\widehat F$ is the minimal extension of the operator $F_0$ (definition \eqref{b12}). Thus, we have:
\begin{equation*}
F_0\subset\widehat F\subset F,
\end{equation*}
i.e. the operator $F$ is an extension of $F_0$.
\end{proof}

A converse statement to Therem \ref{tec1} is valid. We use there the notion
of a locally self-adjoint extension, defined by a Lagrangian submanifold
of the symplectic manifold $(\widehat W,\,\omega)$ (see \eqref{b38}, 
\eqref{b39}, 
Definition \ref{deb2} and Proposition \ref{prb2}).
\begin{theorem}\label{tec2}    
Assume that the conditions (a)-(f) are 
fulfilled, the foliation $\F$ is regular and a locally self-adjoint
extension $F$ of the operator $F_0$ is defined by a Lagrangian 
$C^3$-submanifold $\K_1$ of a symplectic $C^3$-submanifold
$(\widehat W,\;\;\omega)$ of $(\B^2,\,j)$ such that $\widehat W\subset\M$.
Furhermore, assume that $\widehat W$ is $C^3$-diffeomorphic to an open 
domain $\C_1$ of the space $\R^{2m}$ and $\K_1$
is $C^3$-diffeomorphic to a lineal. Then there exists a neighborhood
$U(\K)$ of the submanifold 
$\K = \T(F)$ in $\M$ and a smooth mapping $\Theta$ from
$U(\K)$ into $\R^m$, which satisfy the conditions 
\eqref{c1}-\eqref{c48} of Theorem \ref{tec1}.
\end{theorem}
\begin{proof}
By Theorem \ref{tee1},  there exists a neighborhood
$U(\K_1)$ of the manifold $\K_1$ in 
$\widehat W$ and a mapping
\begin{equation}\label{c49}
\widehat\Theta\in C^1(U(\K_1),\,\R^m),
\end{equation}
\begin{equation}\label{c50}
\K_1 = \left\{\zeta\in \widehat W:\; 
\widehat\Theta(\zeta) = 0\right\}
\end{equation}
and for any point $\zeta\in U(\K_1)$ the derivative
$\widehat\Theta^\prime (\zeta)$ satisfies the conditions:       
\begin{equation}\label{c51}
\I(\widehat\Theta^\prime(\zeta)) = \R^m,
\end{equation}
\begin{equation}\label{c52}
\widehat\Theta^\prime(\zeta)(\Omega(\zeta))^{-1} 
(\widehat\Theta^\prime(\zeta))^* = 0,
\end{equation}
where $\Omega(\zeta)$ is the sqew-self-adjoint operator in
$T_\zeta(\widehat W)$, defined by the 2-form $\omega$ \eqref{b37}, i.e.
\begin{equation*}
\Omega(\zeta) = P_\zeta JP_\zeta
\end{equation*}
and $P_\zeta$ is the orthogonal projection on $T_\zeta(\widehat W)$.
By the property of the manifold $\widehat W$ (Proposition \ref{prb2}) for
any $\zeta\in\widehat W_c$ (Definition \ref{deb1}) there exists an unique
$\widehat\zeta\in\widehat W$, such that
\begin{equation}\label{c53}
\exists\,v\in \B_+\quad
\zeta = G(\widehat\zeta, v).
\end{equation}
Then we can extend the mapping $\widehat\Theta$ on the neighborhood
of $\K = \T(F)$
\begin{equation}\label{c54}
U(\K) = (U(\K_1))_s
\end{equation}
in the following manner:
\begin{equation}\label{c55}
\Theta(\zeta) = \widehat\Theta(\widehat\zeta).
\end{equation}

In view of \eqref{c49} and of $C^1$-smoothness of the mappings $G(\cdot,v)$
(\cite{Zl-3}, Theorem 5.1), the mapping $\Theta$ belongs to the class
$C^1(U(\K),\,R_m)$. Let us take a point
$\widehat\zeta\in U(\K_1)$ and consider the following
subspace of the tangent space $T_{\widehat\zeta}(\widehat W)$:
\begin{equation*}
L_1(\widehat\zeta) = \ker(\widehat\Theta^\prime(\widehat\zeta)).
\end{equation*}
The relations \eqref{c51}, \eqref{c52} and Proposition 4.1 imply that    
$L_1(\widehat\zeta)$  is a Lagrangian subspace of the symplectic 
space $(T_{\widehat\zeta}(\widehat W),\,\omega(\widehat\zeta))$.
Consider the foliation $\F$ defined by \eqref{b24}, and the tangent spaces 
$T_\zeta(\Gamma_\zeta)$ to its leafs $\Gamma_\zeta$. Let us remind
that for each $\zeta\in\M$ the isotropic subspace  
$T_\zeta(\Gamma_\zeta)$ of $(\B^2,\,j)$ is contained in the 
coisotropic subspace $T_\zeta(\M)$ and it is the $j$-
orthogonal complement to $T_\zeta(\M)$ in
$\B^2$. Furthermore, we have by \eqref{b20}:
\begin{equation}\label{c56}
T_\zeta(\Gamma_\zeta) = \T(\D(\E_1(\zeta))).
\end{equation}
Then for any $\widehat\zeta\in U(\K_1)$ the subspace
\begin{equation}\label{c57}
L(\widehat\zeta) = L_1(\widehat\zeta) + T_{\widehat\zeta}
(\Gamma_{\widehat\zeta})
\end{equation}
is Lagrangian in $(\B^2,\,j)$. By the definition of the mapping $\Theta$
(see \eqref{c53}, \eqref{c55}), it is constant along the leafs $\Gamma_\zeta$ 
of the
foliation $\F$. Since $\K = (\K_{1})_s$, then
\eqref{c1} is valid. Furthermore we have:
\begin{equation}\label{c58}
\forall\, \zeta\in U(\K)\quad
T_\zeta(\Gamma_\zeta)\subset\ker(\Theta^\prime(\zeta))
\end{equation}
(see \eqref{c54}), hence we conclude from \eqref{c57}, that for any
$\widehat\zeta\in U(\K_1)$
\begin{equation*}
\ker(\Theta^\prime(\widehat\zeta)) = L(\widehat\zeta).
\end{equation*}
Consequently for any point $\zeta\in U(\K)$, represented 
by the formula \eqref{c53}, we have:
\begin{equation}\label{c59}
\ker(\Theta^\prime(\zeta)) = L(\zeta),
\end{equation}
where
\begin{equation*}
L(\zeta) = G_\zeta(L(\widehat\zeta,\,v).
\end{equation*}
The last equality and symplectomorphic property of the mappings $G(\cdot,v)$
(Proposition \ref{prc1}) imply, that for each $\zeta\in U(\K)$
the subspace $L(\zeta)$ is Lagrangian in $(\B^2,\,j)$, i.e. it is
the graph of a self-adjoint operator $F^\prime(\E_1(\zeta)$.
From the equalities \eqref{c56}, \eqref{c59} and the inclusion 
\eqref{c58} we 
conclude, that $F^\prime(\E_1(\zeta)$ is an extension of
of the symmetric operator $\D(\E_1(\zeta))$, 
which have the defect $m$ (Proposition \ref{prc2} and Definition 
\ref{dec1}).
The last fact, the equality \eqref{c59} and Lemma \ref{lec1} 
imply the conditions
\eqref{c45})-\eqref{c48} of Theorem \ref{tec1}.
\end{proof}

\appendix

\section{Some statements from symplectic
	differential geometry.}
\setcounter{equation}{0}

We  use  in
our  paper  a  characterization  of  a Lagrangian  submanifold  $\K$
of a symplectic manifold  $(\C, \omega)\;(dim(\C) = \ 2m)$ in
terms of a smooth mapping  $\Theta:\;\;\C \rightarrow \R^m$,
whose set of zeros coincides with  $\K$.  In order to establish
this characterization, we need some statements on
extensions of symplectomorphisms.

\subsection{Linear case}

${\bf A.1^o}$. In this subsection we consider
the problem mentioned above
for the linear case. Let $\B^{2m}$ be $2m$-dimensional
real Hilbert space with an inner
product $(\cdot, \cdot)\;(m < \infty)$
and $\Omega$  be a linear invertible skew-self-adjoint
operator in  $\B:\;\;\Omega^\star= - \Omega$. We set
for any $u,v \in\B^{2m}$:
\begin{equation*}
\omega(u,v)=(\Omega u, v).
\end{equation*}
Then the pair $(\B, \omega)$  form a symplectic space.

We will use some well known elementary
notions and facts from
the symplectic linear algebra:
isotropic, coisotropic and Lagrangian subspaces and existence of them;
symplectic basis and existence of
it and other ones ([ ], Chapter  ).
By  $\K^{\perp_\omega}$ we denote the $\Omega$-orthogonal
complement to a subspace
$\K \in \B^{2m}$  in $\B$,
\begin{equation*}
\K^{\perp_\omega} = (\Omega(\K))^\perp.
\end{equation*}
Assume that a subspace  $\K\subseteq\B$ is defined
in the following manner:
\begin{equation}\label{e1}
\K = \ker(\Theta),
\end{equation}
where
\begin{equation}\label{e2}
\Theta \in L(\B^{2m}, \R^{m}).
\end{equation}
Our purpose to obtain a criterion ensuring that the subspace
$\K$  is
Lagrangian. Before we state some auxiliary statements.

\begin{lemma}\label{lme1}
	The following relation is valid:
	\begin{equation}\label{e3}
	\K^{\perp_\omega} = \I(\Omega^{-1}\Theta^\star).
	\end{equation}
\end{lemma}
\begin{proof}
 We have:
\begin{equation*}
\Omega(\K) = \ker(\Theta\Omega^{-1}).
\end{equation*}
Then
\begin{equation*}
\K^{\perp_\omega} = \I((\Theta\Omega^{-1})^\star) =
\I(\Omega^{-1}\Theta^\star).
\end{equation*}
\end{proof}

The following statement is the
straightforward consequence of the
previous lemma:
\begin{lemma}\label{lme2}
	The subspace
	$\K$ (see \eqref{e1}, \eqref{e2}) is coisotropic if and
	only if the following identity holds:
	\begin{equation}\label{e4}
	\Theta \Omega^{-1}\Theta^\star = 0.
	\end{equation}
\end{lemma}

We turn now to a consideration
of Lagrangian subspaces. The
following statement is valid:

\begin{lemma}\label{lme3}
	A subspace $\K\subset\B^{2m}$
	is Lagrangian in  $(\B^{2m},\omega)$, if
	and only if it is coisotropic (isotropic) there and
	\begin{equation}\label{e5}
	\dim(\K) = m.
	\end{equation}
\end{lemma}
\begin{proof}
Let  $P$  be the orthogonal
projector on $\K$. Setting in
\eqref{e3} $\Theta = I - P$, we obtain:
\begin{equation*}
\K^{\perp_\omega}= \I(\Omega^{-1}(I - P)).
\end{equation*}
Then  $\K$  is coisotropic
(isotropic) if and only if the inclusion is
valid:
\begin{equation}\label{e6}
\I(\Omega^{-1}(I - P))\subseteq\I(P)\;\;
(\I(P) \subseteq\I(\Omega^{-1}(I - P))).
\end{equation}
On the other hand, \eqref{e5} is equivalent to:
\begin{equation}\label{e7}
\dim(\I(\Omega^{-1}(I - P))) = \dim(\I(P)) = m.
\end{equation}
The inclusion \eqref{e6} together
with the equality \eqref{e7} are equivalent
to the equality
\begin{equation*}
\I(\Omega^{-1}(I - P)) = \I(P),
\end{equation*}
which means that $\K^{\perp_\omega}= \K$,
i.e. the subspace  $\K$  is Lagrangian. 
\end{proof}

As a consequence
of  Lemma \ref{lme2} and Lemma \ref{lme3} we obtain the
main statement of this subsection:
\begin{proposition}\label{pre1}
	The subspace
	$\K\subseteq\B^{2m}$ (see \eqref{e1}), \eqref{e2}) is
	Lagrangian in $(\B^{2m} \; ,\omega)$  if and only
	if the operator $\Theta$  satisfies
	the condition \eqref{e4} and the relation:
	\begin{equation}\label{e8}
	\I(\Theta)=\R^{m}
	\end{equation}
\end{proposition}
\begin{proof}
The condition \eqref{e8} is
equivalent to \eqref{e5}. Then by
Lemmas \ref{e2}, \ref{e3} we obtain our statement.
\end{proof}

\subsection{Symplectic connection}

 We need a so called symplectic
connection,  which  permit  to  carry  out  "parallel  translations"
of  vectors   on   symplectic
manifolds in similar manner, as
on Riemannian manifolds.  Notice that symplectic connections were studied in \cite{Vais}.

It  will  be convenient for us to consider
a  representation  of  a  symplectic
$C^r$ -smooth manifold ($r\ge 2$) with the
help of charts. Let $\J$  be an open subset of $\R^{m}$,
which  can  be  considered  to  be  a  symplectic
manifold, if a closed non-degenerate
differential 2-form is defined on it:
\begin{equation}\label{e9}
\forall x\in\J, \;\; \forall u, v \in T(\J)\;(\equiv \R^{m})\;\;
\omega(x)(u, v) = (\Omega(x)u, v),
\end{equation}
where
\begin{equation*}
\Omega (x) )\in C^{r-1}(\J, L(\R^{2m})).
\end{equation*}
and for each $x\in\J\;\;\Omega(x)$
is  an  invertible  skew-self-adjoint
operator. For any fixed $x\in\J$
we  consider  the  operator  $\tilde\Omega(x)$,
which is conjugate to the
identity operator $I$  with respect to the
form $\omega(x)$. In other words
$\Omega(x)$  acts from  $\R^{2m}$
into  $(\R^{2m})^\star$   and
\begin{equation}\label{e10}
\forall u, v\in \R^{2m}\;\;
(\tilde\Omega(x)u, v) = \omega(x)(v, u).
\end{equation}
It is easily to see that each
$\tilde\Omega(x)$ is invertible and
\begin{equation*}
\tilde\Omega(x)\in C^{r-1}(\J, L(\R^{2m}, (R^{2m})^\star)).
\end{equation*}
If we identify  $(\R^{2m})^\star$
with $\R^{2m}$,  then $\Omega(x) =(\Omega(x))^\star = - \Omega(x)$,
Consider the following
$C^{r-2}$-mapping from $\J$ into  the  space
$L (\R^{2m}, \R^{2m})$
of antisymmetric  2-forms,  taking  their  values  in
$\R^{2m}$:
\begin{equation}\label{e11}
\forall x \in \J,\;\;  \forall h, \xi \in \R^{2m}\;\;
\Gamma(x)(h, \xi) = - (\tilde\Omega(x))^{-1}(d_x(\Omega(x)(\xi, h))).
\end{equation}
It  turns  out  that  this
mapping  has  the   property   of
a connection, defined by the   2-form
$\omega$, if the last one is
considered to generate a pseudo-Riemannian
metric on $\J$.  Consider  on
$\J$ a $C^r$ -smooth  curve
$\gamma:\;\;x = x(t),\;\; t \in [a,  b]$   and  the
following linear differential
equation of a "parallel  translation"
of vectors  along $\gamma$:
\begin{equation}\label{e12}
\frac{d\xi}{dt} = \Gamma(x(t))(\frac{dx}{dt}, \xi),
\end{equation}
in which  $\xi(t) \in T_x(t)(\J)\;(\equiv \R^{2m})$
for each $t \in [a, b]$.
\begin{proposition}\label{pre2}
	The  evolution
	operator  $U_\gamma(t, t_0)$  of the equation \eqref{e12} realizes
	for  each  $t_0, t\in [a,  b]$   a  linear
	symplectic iomorphism between the
	symplectic  spaces  $(T_x(t_0)(\J),\;\omega(x(t_0)))$  and
	$(T_x(t )(\J),\; \omega(x(t)))$.
\end{proposition}
\begin{proof}
Let us take two arbitrary vectors
$v_1,  v_2\in T_{x(t_0)}(\J)$
and the corresponding solutions of \eqref{e12}:
\begin{equation*}
\xi_k (t) = U_\gamma (t, t_0)v_k\;\;(k = 1, 2).
\end{equation*}
In order to prove our statement,
we ought to show that
\begin{equation}\label{e13}
\frac{d}{dt}[\omega(x(t))(\xi_1 (t), \xi_2 (t))] \equiv 0.
\end{equation}
We have:
\begin{eqnarray}\label{e14}
&&\frac{d}{dt}[\omega(x(t))(\xi_1(t),
	\xi_2(t))] =
	\omega(x(t))(\frac{d\xi_1(t)}{dt}, \xi_2 (t)) +
\nonumber\\
&&  \omega(x(t))(\xi_1 (t),\frac{d\xi_2(t)}{dt})+
\omega_x(x(t))(\xi_1(t),\xi_2 (t)))\cdot\frac{dx(t)}{dt}.
\end{eqnarray}
Taking into account \eqref{e10}, \eqref{e11}
and definition \eqref{e9} of the
operator  $\tilde\Omega$, we have:
\begin{equation}\label{e15}
\omega(x(t))(\frac{d\xi_1(t)}{dt},\,\xi_2(t)) =
- \omega_x (x(t))(\xi_1 (t), {dx(t)}{dt}))\cdot\xi_2(t),
\end{equation}
\begin{equation}\label{e16}
\omega(x(t))(\xi_1 (t),\, \frac{d\xi_2(t)}{dt}) =
\omega_x (x(t))(\xi_2 (t), \frac{dx(t)}{dt})\cdot\xi_1 (t).
\end{equation}
From the identities \eqref{e14} - \eqref{e16} we
obtain, taking into  account
that $\omega$  is a closed  2-form:
\begin{eqnarray*}
&&\frac{d}{dt}[\omega(x(t))(\xi (t), \xi (t))] =
	\omega_x (x(t))(\xi_2 (t), \frac{dx(t)}{dt}))\cdot\xi_1 (t) -
\nonumber\\
&&\omega_x (x(t))(\xi_1 (t), \frac{dx(t)}{dt}))\cdot\xi_2(t) +
\omega_x (x(t))(\xi_1 (t), \xi_2 (t)))\cdot\frac{dx(t}{dt} =
\nonumber\\
&& d\omega(x(t))\big(\xi_1 (t), \xi_2 (t), \frac{dx(t)}{dt}\big) = 0,
\end{eqnarray*}
i.e. the identity \eqref{e13} is valid.
\end{proof}

\subsection{Extension of a diffeomorphism to a tubular neighborhood}

We need some statements on an extension
of  a  diffeomorphism
between two manifolds to a diffeomorphism
between  their  tubular neighborhoods. Before we state a
statement on a  trivialization  of some class of vector bundles.

\begin{lemma}\label{lme4}
	Let $\{N_l\}_{l\in E}$ be a family of
	subspaces of  a  Hilbert space  $\B$  ($E$ is a Banach  space),
	such  that  the  orthogonal
	projectors  P(l)  on them have the property:
	\begin{equation*}
	P(\cdot) \in C^r(E, L(\B)).
	\end{equation*}
	Then there  exists  a  family
	$\{U(l)\}_{l\in E}$ of  topological  linear
	isomorphisms from $N_0$   onto  $N_1$, such that
	\begin{equation*}
	U(\cdot) \in C^{r-1}(E, L(N_0 , \B)).
	\end{equation*}
	In other words, the mapping
	\begin{equation}\label{e18}
	\U(l, v) = \{l, U(l)v\}\;\; (l \in E, v \in N_0 )
	\end{equation}
	realizes a $C^{r-1}$- smooth  BP-morphism
	between  the  trivial  vector
	bundle $E \times N_0\rightarrow E$  and the vector bundle,
	defined  by  the  family
	$\{N_l\}_{l\in E}$.
\end{lemma}
\begin{proof}
We set $Q(l) = I - P(l)$.
For any fixed $l \in  E$  consider
the following evolution equation in the Hilbert space  $\B$:
\begin{equation*}
\frac{d\xi(t)}{dt} = [\frac{dP(tl)}{dt}P(tl)
+ \frac{dQ(tl)}{dt})Q(tl)]\xi.
\end{equation*}
It is known  that  the
evolution  operator $V_l(t, t_0)$  of  this
equation is unitary and it
"turns" the subspaces $N_l$ and $N_l^\perp$ , i.e.,
\begin{equation*}
P(tl)V_l(t, t_0 ) = V_l(t, t_0)P(t_0l),\;\;
Q(tl)V_l(t, t_0 ) = V_l(t, t_0)Q(t_0l)
\end{equation*}
(see [ ], Chapt.    ). It is easily to
see that the mapping
\begin{equation*}
U(l) = V_l(1, 0))\vert_{N_0}
\end{equation*}
has desired properties.
\end{proof}

The following  statement  on
extension  of  diffeomorphisms  is
valid:
\begin{proposition}\label{pre4}
	Let   $\K$
	be  a  $C^r$- smooth $(r\ge 3)$  regular
	submanifold of a Hilbert  space $\B$.
	Assume  that  a  mapping  $U_\K$  realizes a $C^r$- diffeomorphism
	between a subspace $\hat\K$  of  $\B$  and  the
	submanifold  $\K$. Then it can be
	extended to a $C^{r-2}$-diffeomorphism
	$U_\X$   between a tubular
	neighborhood  $\X(\K)$   of   $\K$  and  a  tubular
	neighborhood  $\X(\hat\K)$  of $\hat\K$.
\end{proposition}
\begin{proof}
Consider the normal
bundle $\{N_x\}_{x\in\K}$ on  the  submanifold
$\K$:
\begin{equation*}
\forall x \in \K:\;\;  N_x  = (T_x (\K))^\perp.
\end{equation*}
Furthermore, consider  on $\B$   the
geodesic  flow,  generated  by
$\B$-metric and the exponential
mapping $Exp$  defined  by  this  flow
\cite{Lng}. Let  $\X(\K)$
be the tubular  neighborhood  defined
by the normal bundle  $N(\K)$
and the mapping $Exp$, i.e.,
\begin{equation*}
T(\K) = f_N(Z_N ),
\end{equation*}
where
\begin{equation*}
f_N  = Exp\vert_{N(\K)},
\end{equation*}
$Z_N$ is a neghborhood of the zero
section  $\zeta_\K$ of $N(\K)$ and $f_N$
realizes a $C^{r-2}$-diffeomorphism between
$Z_N$ and  $\X(\K)$ . Let us transfer the bundle  $N(\K)$
from the base $\K$  to the base $\hat\K$  with the help of 
mapping $U_\K$ , i.e. consider the  following
vector bundle:
\begin{equation}\label{en}
\hat N(\hat\K) = (U_\K^{-1})^\star(N(\K)).
\footnote{Here and in the sequel $f^\star$ denotes the lifting of a vector
	bundle $E\rightarrow X$ from the base $X$ to
	the base $\hat X$ with the help of a mapping $f:\;\;\hat X\rightarrow X$
\cite{Lng}. We hope that
	this notation will not cause
	a confusion with the notation of adjoint linear
	homomorphisms.}
\end{equation}
Let us notice that $U_\K^\star$ is a BP-isomorphism
and $U_\K^{-1})^\star)^{-1}= U_\K^\star$.
By  Lemma  \ref{lme4}  there  exists  a $C^{r-1}$-smooth
BP-isomorphism  $\U$  between
the trivial vector bundle  $\hat\K\times N_0\rightarrow\hat\K$
and  the vector bundle $N(\K)$.
Since  ${\mathrm codim}(\hat\K) ={\mathrm codim}(\K)$,  then  there
exists a topological linear
isomorphism $S$ between $\tilde N = (\hat\K)^\perp$   and $N_0$.
Let $P_{\tilde N}$  be the orthogonal
projection on $\tilde N$. Then the mapping
\begin{equation*}
\tilde U(x) = \{(I - P_{\tilde N})x, S(P_{\tilde N}(x))\}
\end{equation*}
realizes a topological
linear isomorphism between $\B$ and $\hat\K\times N_0$.

Consider the mapping  $\Phi = \tilde\U^{-1}\U^{-1}U_\K^\star f_N^{-1}$,
which  realizes  a $C^{r-2}$-diffeomorphism between $\X(\K)$
and a tubular neighborhood  $\X(\hat\K)$  of $\hat\K$.
Then the mapping $U_\X  = \Phi^{-1}$  is the
desired extension of $U_\K$. 
\end{proof}

We need also the following statement
on a  diffeomorphism  of  a neighborhoods of a subspace
$\K$ of a Hilbert  space $\B$,  generated
by a BP-automorphism of the restriction of the tangent
bundle  $T(\B)$ on $\K$.  Let  $\K^\prime$  
be a complement of  $\K$  in $\B$, 
$P_\K$   be the projection on  $\K$  along to  $\K^\prime$ ,  
$Q_\K  = I - P_{\K^{\prime}}$.
\begin{proposition}\label{pre5}
	Let  $W$  be  a  $C^r$ -smooth  BP-automorphism  of        
	the vector bundle  $T(\B)\vert_\K$, such that
	\begin{equation*}
	W\vert_{T(\K)} = id_{T(\K)}.
	\end{equation*}
	Then there exists a $C^r$-diffeomorphism $\V$ of a neighborhood  
	$\X_1(\K)$ of $\K$ onto a neighborhood $\X_2(\K)$ of $\K$, 
	such that
	\begin{equation}\label{f1}
	\V\vert_\K  = id_\K             
	\end{equation}
	and
	\begin{equation}\label{f2}
	T(\V)\vert_{T(\B)\vert_\K} = W.
	\end{equation}    
\end{proposition}
\begin{proof}
For any $x \in  \K$ 
the  BP-automorphism  $W$ has the following form on  the        
fiber $T_x(\B)\;(\equiv \B)$:
\begin{equation*}
W(x) = Id_\K  + U(x),
\end{equation*}
where 
\begin{equation*}
U(x) = W(x)\vert_{\K^\prime}.
\end{equation*}
Furthermore, we have:
\begin{equation}\label{f3}
W(x) \in Laut(\B).                   
\end{equation}
Consider the following mapping of $\B$ into itself:
\begin{equation*}
\tilde\V(x) = P_\K x + U(P_\K x)Q_\K x.
\end{equation*}
It is obviously that \eqref{f2} is valid.  Let  us  calculate  the 
derivative of $\V(x)$ along any direction $\xi \in \B$:
\begin{equation*}
\V^\prime(x)\xi = P_\K \xi + U_x(P_\K x)Q_\K x)\cdot P_\K \xi + 
U(P_\K x)Q_\K \xi.
\end{equation*}
In particular, for  $x\in\K$:
\begin{equation*}
\V^\prime(x)\xi = P_\K\xi + U(P_\K x)Q_\K\xi = W(x)\xi,
\end{equation*}
hence \eqref{f3} holds. Furthermore, 
in view of \eqref{f4}, the mapping  
$\tilde\V$ is a local 
$C^r$-diffeomorphism at each point $x\in\K$. 
There is only to 
construct a tubular neighborhood 
$\X_1(\K)$  of $\K$, such that the 
restriction
\begin{equation*}
\V = \V\vert_{\X_1(\K)}
\end{equation*}
realizes a $C^r$ -diffeomorphism on its range $\X_2(\K)$. 
\end{proof}

\subsection{Main statement}

In this subsection we 
shall state the main statement  of  this 
section, which give a description of a Lagrangian manifold in terms 
of a mapping, whose set of zeros coincides with  it.  We  eslablish 
before a statement on an extension of a  symplectomorphism  from  a 
Lagrangian manifold to its tubular neighborhood. This statement can 
be considered to be a global variant of Givental Theorem  for  the 
case of Lagrangian manifolds (\cite{Tr}, Chapt VII, Sect. 3).           

\begin{proposition}\label{pre6}
	Let $(\C,\omega)$ be a symplectic $C^r$ -manifold 
	($r\ge 3$) 
	of a finite dimension  $2m$ and $\K$ be its Lagrangian 
	submanifold. Assume that $\C$ is $C^r$ -diffeomorphic 
	to an open  domain  
	$\C$ of the space $\R^{2m}$ with the canonical 
	symplectic form $j$ and $\K$  
	is $C^r$ -diffeomorfhic to a 
	Lagrangian subspace $\hat\K$ of $(\R^{2m},\,j)$. Then 
	there exists a $C^{r-2}$-diffeomorphism  $\Psi$ 
	from a tubular  neighborhood 
	$\X(\hat\K)$ of $\hat \K$ onto a tubular 
	neighborhood $\X(\K)$ of $\K$, such that:
	\begin{equation}\label{f4}
	\Psi(\hat\K) = \K
	\end{equation}
	and
	\begin{equation}\label{f5}  
	\Psi^\star(\omega\vert_{T(\X(\K))}) = j\vert_{T(\X(\hat\K))}.
	\end{equation}
\end{proposition}
\begin{proof}
Assume that we have constructed a $C^{r-1}$ -diffeomorphism 
$\Y$ 
from a tubular neighborhood $\X_0(\hat\K)$ of $\hat\K$  
onto a tubular 
neighborhood $\X_0(\K)$ of $\K$, such that
\begin{equation}\label{f6}
\Y(\hat\K) = \K
\end{equation}
and
\begin{equation}\label{f7}
\Y^\star(\omega\vert_{T(\X(\K))\vert_\K}) = 
j\vert_{T(\X(\hat\K))\vert_{\hat\K}}.
\end{equation}  
Then setting 
\begin{equation*}
\tilde\omega =\Y^\star(\omega\vert_{T(\X(\K))}),
\end{equation*}
we have by \eqref{f7} that
\begin{equation*}
\tilde\omega\vert_{T(\X(\K))\vert_\K} = 
j\vert_{T(\X(\hat\K))\vert_{\hat\K}}.
\end{equation*}
Then by the Veinstein theorem ([ ], Chapt.          )  there  exist 
tubular neighborhoods $\X(\hat\K)\subseteq\X_0(\K),\;\;
\X_1(\hat\K)\subseteq\X_0(\hat\K)$ of $\hat\K$ and a 
$C^{r-1}$-diffeomorphism $\Phi$ from $\X(\hat\K)$ 
onto $\X_1(\K)$, such that
\begin{equation*}
\Phi\vert_{\hat\K}= id_{\hat\K},\;\;
\Phi\star(\omega\vert_{T(\X_1(\hat\K))}) = 
j\vert_{T(\X(\hat\K))}.
\end{equation*}
Then the mapping  $\Psi = \Y\circ\Phi$  satisfies 
the conditions \eqref{f4},  \eqref{f5} 
with $\X(\K) = \Psi(\X(\hat\K))$.

We turn now to the construction of the  mapping  $\Y$ satisfying 
the conditions \eqref{f6}, \eqref{f7}. 
Before let us map the manifold $\C$  
onto the open domain $\tilde\C$ of $\R^{2m}$ by means of a 
$C^r$ -diffeomorphism  
$U_\C$ . We set $\tilde\K = U_\C(\K)$. 
Let $U_\K$ be a $C^r$ -diffeomorphism from $\tilde\K$  
onto the Lagrangian subspace $\hat\K \subseteq  \R^{2m}$. 
Then by Proposition \ref{pre4} 
there exists a $C^{r-2}$-diffeomorphic extension  
$U_\X$ of $(U_\K^){-1}$, which maps a tubular neighborhood        
$\X_2(\hat\K)$ of $\hat\K$ onto a 
tubular neighborhood $\X_3(\tilde\K)$ of $\tilde\K$. 
We set 
\begin{equation}\label{f8}
\X_3 (\K) = (U_\C)^{-1}(\X_3(\hat\K)),\nonumber\\
\hat\omega = (U_\C^{-1}\circ U_\X)\star(\omega\vert_{T(\X_3(\K))}).
\end{equation}
Since $U_\C^{-1}\circ U_\X(\hat\K) = \K$  
and $\K$ is Lagrangian in $(\X_3(\K), \omega)$, then $\hat\K$  
is Lagrangian in $(\X_2(\hat\K), \hat\omega)$.

Consider on the symplectic manifold 
$(\X_2(\hat\K), \hat\omega)$  the  symplectic 
connection  $\Gamma(x)(\cdot, \cdot)$ 
(see \eqref{e11}). We shall translate vectors 
with the help of this connection along the rays
\begin{equation*}
\gamma_l:\;\; x =tl,\;t \ge  0,\; l \in \K.
\end{equation*}
In other words, we consider the family of differential equations:
\begin{equation}\label{f9}
\frac{d\xi}{dt} = \Gamma(tl))(l, \xi).
\end{equation}
For a fixed $l \in \K$ we denote by $U_l(t, t_0)$ 
the evolution operator 
of the equation \eqref{f9}. Let us show that $\hat\K$ is 
invariant  with 
respect of this operator. This is equivalent to  the  fact  that  for 
each fixed $l\in\K,\;\;t \ge 0$ the 
restriction  of  the  vector  field        
$\{\Gamma(tl))(l, \xi)\}_{\xi\in\X_2(\hat\K)}$  on $\hat\K$ 
is a vector field on $\hat\K$. This means 
that 
\begin{equation}\label{f10}
\forall \xi \in\hat\K:\;\;  
\Gamma(tl))(l, \xi)\in T(\hat\K)\;(\equiv\hat\K).
\end{equation}
Let us prove the last property. Since $\hat\K$  is Lagrangian,  then  for 
any $l \in\hat\K$  and $\xi, h \in T_{tl}(\hat\K)$
\begin{equation}\label{f11}
\omega(tl)(h, \xi) = 0,    
\end{equation}
in particular,
\begin{equation*}
\omega(tl)(l, \xi) = 0.
\end{equation*}
The last identity implies that ,
\begin{equation}\label{f12}
\forall h \in T_{tl}(\hat\K):\;\;  
d_l(\omega(tl)(l, \xi))\cdot h =\nonumber\\
= td_x\omega(tl)(l, \xi)\cdot h + \omega(tl)(h, \xi) = 0. 
\end{equation}
In view of \eqref{f11}, the identity \eqref{f12} 
is equivalent to the following one:
\begin{equation}\label{f13}
\forall h \in T_{tl}(\hat\K):\;\;d_x\omega(tl)(l, \xi)\cdot h = 0.
\end{equation}
According to the definition of the operator 
$\tilde\Omega(x)$ \eqref{e10}   and 
of $\Gamma(x)(\cdot, \cdot)$ \eqref{e11} 
the relation $y = \Gamma(tl))(l, \xi)$  is equivalent        
to
\begin{equation}\label{f14}
\forall h \in T_{tl}(\R^{2m})\;(\equiv \R^{2m}):\nonumber\\
\omega(tl)(y, h) = - d_x\omega(tl)(l, \xi)\cdot h.
\end{equation}
Then \eqref{f13}, \eqref{f14} and the fact that $\hat\K$ 
is a Lagrangian 
submanifold of 
the symplectic manifold $(\X_2(\hat\K), \hat\omega)$,  
imply the inclusion: $y \in T_{tl}(\hat\K)$.  So, 
the invariance of $\hat\K$ with 
respect to the evolution operator  of 
the equation \eqref{f9} is proved.

We now turn to a symplectomorphic trivialization of the vector 
bundle $T(\X_2(\hat\K))\vert_{\hat\K}$ 
by means of the above symplectic connection. Let        
us set
\begin{equation*}
\W(\xi_0, l) = \{l, W_l \xi_0\}, 
\end{equation*}
where
\begin{equation*}
W_l  = U_l(1, 0).
\end{equation*}
It is easy to see that the mapping $\W$ realizes 
a $C^{r-2}$-smooth  
BP-morphism of the trivial vector bundle  
$\hat\K\times T_0(\R^{2m})\rightarrow\hat\K$  onto the 
vector bundle $T(\X_2(\K))\vert_{\hat\K}$. 
The  invariance  of  $\hat\K$, 
which  have 
been proved above, implies the following property of  $\W$:
\begin{equation}\label{f15}
\W(\hat\K \times T_0(\K)) = T(\hat\K).
\end{equation}
Furthermore,  by  the  property  of  the  symplectic   connection 
(Proposition \ref{pre2}) on each fibre $\{l\}\times T_0(\R^{2m})$ 
the  mapping $\W$  
realizes a linear symplectic isomorphism between the symplectic 
spaces $(\{l\} \times T_0(\R^{2m}), \hat\omega(0))$ 
and \\$(T_l(\X_2(\hat\K)), \hat\omega(l))$.
For a convenience we shall identify the spaces  
$T_l(\R^{2m})$ with  
$\R^{2m}$. We choose a symplectic basis in the symplectic space
$(\R^{2m}, j)$ (respectively, in  $(\R^{2m}, \hat\omega(0))$) 
in the following manner:
\begin{equation*}
e_1^{(0)}, e_2^{(0)},\dots, e_m^{(0)}, 
f_1^{(0)}, f_2^{(0)}, \dots, f_m^{(0)}
\end{equation*}
and, respectively,
\begin{equation*}
e_1^{(0)}, e_2^{(0)},\dots, e_m^{(0)},
\hat f_1^{(0)},\hat f_2^{(0)}, \dots,\hat f_m^{(0)},
\end{equation*}
where $e_k\in\hat\K\;(k = 1,2,  ...,  m)$.  Let  us  remind  that   the 
subspace $\hat\K$ is Lagrangian 
in the both symplectic spaces.  Consider 
the symplectic linear  diffeomorphism $S$ 
from $(\R^{2m},  j)$ onto        
$(\R^{2m}, \hat\omega(0))$, 
defined by the correspondence of the above basises. It is 
clear that
\begin{equation}\label{f16}
S\vert_{\hat\K} = id_{\hat\K}.             
\end{equation}
Then for each $l\in\hat\K$ the mapping $W_l\circ S$ 
realizes a symlectomorphism between 
$(\R^{2m},  j)$ and $(T_l(\X_2(\hat\K)), \hat\omega(l))$.

Our aim is to construct for any $l \in \K$ a symplectomorphism  
$\tilde W_l$ between the above symplectic spaces, such that the mapping
\begin{equation}\label{f17}
\tilde\W(l, v) = {l,\tilde W_l v}
\end{equation}
realizes a $C^{r-2}$-smooth BP-morphism between the vector bundles
\begin{equation*}
T(\R^{2m})\vert_{\hat\K},\;\;T(\X_2(\hat\K))\vert_{\hat\K}
\end{equation*}
and, moreover, the relation holds
\begin{equation}\label{f18}
\W\vert_{T(\hat\K))}  = id_{T(\hat\K))},      
\end{equation}
if we identify on the both 
manifolds $\R^{2m},\;\;\X_2(\hat\K)$ the vector 
bundles $T(\hat\K)$ and $\hat\K \times\hat \K$. 
Taking into account the last remark,  we 
can consider $e_1^{(0)}, e_2^{(0)},\dots, e_m^{(0)}$ 
to be a basis in each tangent        
space $T_l(\hat\K)$. In view of \eqref{f15}, \eqref{f16},
\begin{equation*}
\forall l \in \hat\K:\;\; W_l \circ S(T_l(\hat\K)) = T_l(\K).
\end{equation*}
Then the vector functions
\begin{equation*}
e_k (l) = (W_l\circ S)^{-1}(e_k^{(0)})\; (k = 1, 2, ..., m)  
\end{equation*}
are $C^{r-2}$-smooth and for each fixed $l\in\hat\K$ 
they form a basis in the tangent 
space $T_l(\hat\K)$ identified with $\hat\K$. 
Consider in  $\R^{2m},\; j)$  the 
following Lagrangian subspace $\hat\K^\prime$, 
which is complementary to $\hat\K$:
\begin{equation*}
\K  = {\rm span}(f_1^{(0)}, f_2^{(0)}, \dots, f_m^{(0)}).  
\end{equation*}
Let us construct a basis in $\hat\K^\prime$ for each fixed 
$l \in\hat \K$
\begin{equation*}
\tilde f_1(l), \tilde f_2(l), \dots, \tilde f_m(l),
\end{equation*}
such that the vectors
\begin{equation*}
e_1 (l), e_2 (l), \dots, e_m (l), 
\tilde f_1 (l),\tilde f_2 (l), \dots,\tilde f_m (l)
\end{equation*}
form a symplectic basis in $(\R^{2m}   , j)$. 
The vectors $\tilde f_i(l)$ have the form:
\begin{equation*}
\sum_{k=1}^ma_{ik}f_k^{(0)},
\end{equation*}
where for each fixed $i$ the numbers $a_{ik}\;(k =1, 2, ..., m)$  form 
the solution of the following linear system:
\begin{equation*}
\sum_{k=1}^m a_{ik}j(f_k^0, e_s (l)) = \delta_{is}\;(s = 1, 2, ..., m).
\end{equation*}
It is easy to see that this system has the  unique  solution  and 
the vector functions  $\tilde f_i (l)$  are $C^{r-2}$-smooth. 

Consider the sequence of vectors:
\begin{equation}\label{f19}
e_1^{(0)}, e_2^{(0)},\dots, e_m^{(0)},
f_1(l),f_2(l),\dots,f_m(l),
\end{equation}
where
\begin{equation*}
f_i(l) = W_l\circ S(\tilde f_i (l)). 
\end{equation*}
Since  $e_i^{(0)} = W_l\circ S(e_i(l))$  and $W_l\circ S$  
realizes a linear  symplectic 
isomorphism between $(\R^{2m}, j)$  and $(T_l(\R^{2m}),\hat\omega(l))$, 
then the system \eqref{f19} form a symplectic basis in the last 
symplectic space. Let us 
construct the linear symplectic isomorphism $\hat W_l$ 
between $(\R^{2m},  j)$        
and $(T_l(\R^{2m}),\hat \omega(l))$  by 
means of the correspondence  of  the  symplectic bases
\begin{equation*}
\tilde W_l(e_i^{(0)}) = e_i^{(0)},\;\; 
\tilde W_l(f_i^{(0)}) = f_i(l)\;\; (i = 1, 2, ..., m).
\end{equation*}
Then the mapping of the form \eqref{f17} 
defines a $C^{r-2}$-smooth 
BP-morphism $\W$ between the vector  bundles 
$T(\R^{2m})\vert_{\hat\K}$ and 
$T(\X_2(\hat\K))\vert_{\hat\K}$ having the property 
\eqref{f18}. Furthermore, this 
isomorphism is symplectic, i.e.,
\begin{equation}\label{f20}
\tilde\W^\star(\omega\vert_{T(\X_2(\hat\K)))\vert_{\hat\K}}) = 
j\vert_{T(\R^{2m})\vert_{\hat\K}}.   
\end{equation}
By Proposition \ref{pre5} we can construct a 
diffeomorphism $\V$  from a 
tubular neighborhood $\X(\hat\K)$ of 
$\hat\K$  onto a tubular  neighborhood 
$\X_1(\tilde\K) \subseteq \X_3(\tilde\K)$ of $\tilde\K$, 
such that
\begin{equation*}
\V\vert_{\hat\K}  = id_{\hat\K} 
\end{equation*}
and
\begin{equation*}
T(\V)\vert_{T(\X_1(\tilde\K)))\vert_{\tilde\K}} = \tilde\W.
\end{equation*}
The last relations and the equalities \eqref{f8}, \eqref{f20} 
imply that the 
mapping $\Y = U_\C\circ U_\X\circ\V$ has the properties 
\eqref{f6}, \eqref{f7}. 
\end{proof}

Let us turn now to the main statement of this section.

\begin{theorem}\label{tee1}
	Let $(\C, \omega)$ be a symplectic $C^r$ -manifold of a 
	finite dimension $2m$ endowed by a Riemannian metric 
	$\{(\cdot, \cdot)_x\}_{x\in\C}$.  A 
	non-empty subset $\K \subseteq \C$ 
	is a Lagrangian submanifold of $(\C, \omega)$, if there exists a tubular neighborhood $\X(\K)$ of $\K$ 
	and a mapping $\Theta$ 
	from $\X(\K)$ into the space $\R^{m}$ satisfying the conditions:
	\begin{equation}\label{f21}
	\Theta \in C^r(\X(\K), \R^{m}), 
	\end{equation}
	\begin{equation}\label{f22}
	\K = \{x \in \X(\K):\;\;\Theta(x) = 0\}, 
	\end{equation}
	\begin{equation}\label{f23}
	\forall x \in \X(\K):\;\;\I(\Theta^\prime)(x)) = \R^{m}, 
	\end{equation}
	\begin{equation}\label{f24}
	\forall x \in \X(\K):\;\;\Theta^\prime(x)(\Omega(x))^{-1}
	(\Theta^\prime(x))^\star  = 0,  
	\end{equation}
	where the operator $(\Theta^\prime(x))^\star$ 
	is conjugate to the operator $\Theta^\prime(x)$ 
	(if the latter is considered to be acting from the Hilbert space        
	$(T_x(\C), (\cdot, \cdot)_x)$ into the Hilbert space 
	$\R^{m}$) and the 
	operator $\Omega(x)$ is a skew-self-adjoint invertible  operator 
	generating in the first space the 2-form $\omega(x)(\cdot,\cdot)$, i.e.,
	\begin{equation*}
	\forall u, v \in T_x(\C):\;\;\omega(x)(u,v) = (\Omega(x)u, v)_x.
	\end{equation*}
	If $r \ge 3$, the manifold $\C$ is  $C^r$ -diffeomorphic 
	to an open domain  
	$\tilde\C$ of the space $\R^{2m}$ and $\K$ is $C^r$ -diffeomorphic  
	to  a  linear 
	space, then the necessary condition for $\K$ to be a Lagrangian 
	submanifold of $(\C, \omega)$  
	is the existence of a mapping $\Theta$, such that 
	the conditions \eqref{f21} - \eqref{f24}  are  satisfied for it 
	with $r-2$ instead of $r$.
\end{theorem}
\begin{proof}
The first part of our theorem follows from Proposition 
\ref{pre1} applied to the tangent space at each point of $\X(\K)$. 
Let us prove the second part. According to the condition imposed on 
the Lagrangian submanifold $\K$, there  exists  a  $C^r$ -diffeomorphism 
between it and a Lagrangian subspace $\K$ of the symplectic space 
$(\R^{2m}, j)$. Then by Proposition \ref{pre5} 
there exists a 
$C^{r-2}$-diffeomorphism $\Psi$ from a tubular neighborhood 
$\X(\hat\K)$ of $\hat\K$  
onto a tubular neighborhood $\X(\K)$ of $\K$ satisfying the conditions 
\eqref{f4}, \eqref{f5}. Let $P$ be the orthogonal projection  
on the subspace $\hat\K \subseteq \R^{2m}$, i.e.,
\begin{equation}\label{f25}
\hat\K =\ker(I - P).                    
\end{equation}
Since $\hat\K$ is Lagrangian, then by Proposition \ref{pre1}, 
in  which $\R^{m}=\hat\K^\perp$, the conditions are fulfilled:
\begin{equation}\label{f26}
rank(P) = m,                     
\end{equation}
\begin{equation}\label{f27}
(I - P)J(I - P) = 0.                 
\end{equation}
We set $\Theta = (I - P)\Psi^{-1}$. 
Taking into account \eqref{f25} - \eqref{f27}, we 
obtain, that this mapping  $\Theta$  
satisfies the desired conditions.
\end{proof}

\end{document}